 \documentclass[12pt]{article}

 \usepackage{booktabs}
 \usepackage{threeparttable}
 \usepackage{mathrsfs,amsfonts,amsmath}
 \usepackage{harvard}
 \usepackage[dvips]{color}

 \newtheorem{theorem}{Theorem}[section]
 \newtheorem{definition}[theorem]{Definition}
 \newtheorem{lemma}[theorem]{Lemma}
 \newtheorem{corollary}[theorem]{Corollary}
 \newtheorem{proposition}[theorem]{Proposition}
 \newtheorem{remark}[theorem]{Remark}
 \newtheorem{condition}[theorem]{Condition}
 \newtheorem{example}{Example}[section]

 \def\bdefinition{\begin{definition}\rm{}\def\edefinition{\end{definition}}}
 \def\blemma{\begin{lemma}}\def\elemma{\end{lemma}}
 \def\bproposition{\begin{proposition}}\def\eproposition{\end{proposition}}
 \def\ttheorem{\begin{theorem}}\def\etheorem{\end{theorem}}
 \def\bcorollary{\begin{corollary}}\def\ecorollary{\end{corollary}}
 \def\bremark{\begin{remark}}\def\eremark{\end{remark}}
 \def\bcondition{\begin{condition}}\def\econdition{\end{condition}}

 \def\benumerate{\begin{enumerate}}\def\eenumerate{\end{enumerate}}
 \def\bitemize{\begin{itemize}}\def\eitemize{\end{itemize}}
 \def\itm{\item}

 \def\beqlb{\begin{eqnarray}}\def\eeqlb{\end{eqnarray}}
 \def\beqnn{\begin{eqnarray*}}\def\eeqnn{\end{eqnarray*}}
 \def\nnm{\nonumber}\def\ccr{\nnm\\}\def\ar{\!\!\!&}

 \def\mcr{\mathscr}\def\mbb{\mathbb}\def\mbf{\mathbf}\def\mrm{\mathrm}

 \def\proof{\noindent{\it Proof.~~}}\def\qed{\hfill$\Box$\medskip}

 \def\b{\mrm{b}}

\begin{document}

~

\bigskip\bigskip\bigskip

\centerline{\Large\bf Asymptotic properties of estimators in}

\smallskip

\centerline{\Large\bf a stable Cox-Ingersoll-Ross model}

\bigskip

\centerline{By Zenghu Li\footnote{ Supported by NSFC, 973 Program and 985
Program.} and Chunhua Ma\footnote{ Corresponding author. Supported by
NSFC and CSC.}}

\medskip

\centerline{Beijing Normal University and Nankai University}

\bigskip

{\narrower

\noindent\textit{Abstract:} We study the estimation of a stable
Cox-Ingersoll-Ross model, which is a special subcritical continuous-state
branching process with immigration. The process is characterized in terms
of some stochastic equations. The exponential ergodicity and strong mixing
property of the process and the heavy tail behavior of some related random
sequences are studied. We also establish the convergence of some point
processes and partial sums associated with the model. From those results,
we derive the consistency and central limit theorems of the conditional
least squares estimators and the weighted conditional least squares
estimators of the drift parameters based on low frequency observations. A
weakly consistent estimator is also proposed for the volatility
coefficient based on high frequency observations.

\smallskip



~

\par}

\bigskip


\section{Introduction}

\setcounter{equation}{0}

The \textit{Cox-Ingersoll-Ross model} (CIR-model) introduced by Cox et
al.\ (1985) has been used widely in the financial world. This model has
many appealing advantages. In particular, it is mean-reverting and remains
positive. Let $a>0$, $b>0$ and $\sigma>0$ be given constants. The
classical CIR-model is a positive diffusion process $\{X(t): t\ge 0\}$
defined by
 \beqlb\label{1.1}
dX(t) = (a-bX(t))dt + \sigma\sqrt{X(t)} dB(t),
 \eeqlb
where $\{B(t): t\ge 0\}$ is a standard Brownian motion. The process
defined by (\ref{1.1}) has continuous sample paths and light tailed
marginal distributions.

It is well-known that many financial processes exhibit discontinuous
sample paths and heavy tailed distributions. Let $(\Omega,\mcr{F},
\mcr{F}_t,\mbf{P})$ be a filtered probability space satisfying the usual
hypotheses. A natural generalization of (\ref{1.1}) is the stochastic
differential equation
 \beqlb\label{1.2}
dX_t = (a-bX_t)dt + \sigma\sqrt[\alpha]{X_{t-}} dZ_t,
 \eeqlb
where $\{Z_t: t\ge 0\}$ is a spectrally positive stable
$(\mcr{F}_t)$-L\'evy process with index $1<\alpha\le 2$. For $\alpha=2$,
we understand the noise as a standard Brownian motion, so (\ref{1.2})
reduces to (\ref{1.1}). When $1<\alpha<2$, we assume it is a stable
process with L\'{e}vy measure
 \beqlb\label{1.3}
\nu_\alpha(dz) :=
\frac{1_{\{z>0\}}dz}{\alpha\Gamma(-\alpha)z^{\alpha+1}}.
 \eeqlb
By a result of Fu and Li (2010), there is a pathwise unique positive
strong solution $\{X_t: t\ge 0\}$ to (\ref{1.2}). We refer to this process
as a \textit{stable Cox-Ingersoll-Ross model} (SCIR-model). We shall see
that the discontinuous SCIR-model indeed captures the important heavy tail
property. The reader may refer to Borkovec and Kl\"uppelberg (1998),
Embrechts et al.\ (1997, Section~7.6) and Fasen et al.\ (2006) for similar
modifications of the CIR-model. The SCIR-model is a particular form of the
so-called \textit{continuous-state branching processes with immigration}
(CBI-processes), which arise as scaling limits of \textit{Galton-Watson
branching processes with immigration} (GWI processes); see, e.g., Kawazu
and Watanabe (1971). The general CBI-processes were also constructed and
studied in terms of stochastic integral equations in Dawson and Li (2006,
2012), Fu and Li (2010) and Li and Ma (2008).

The estimation for stochastic processes based on the minimization of a sum
of squared deviations about conditional expectations was developed in
Klimko and Nelson (1978). They applied their results to the
\textit{conditional least squares estimators} (CLSEs) of the offspring and
immigration means of subcritical GWI processes. Their estimators are
essentially the same as those studied by Quine (1976, 1977). By the
results of Klimko and Nelson (1978) and Quine (1976, 1977), under a finite
third moment condition, as the sample size $n$ goes to infinity, the
errors of the CLSEs decay at rate $n^{-1/2}$ and they are asymptotically
Gaussian; see also the earlier work of Heyde and Seneta (1972, 1974). The
asymptotic properties of CLSEs of GWI processes with general offspring
laws were studied in Venkataraman (1982) and Wei and Winnicki (1989).
Based on the idea of Nelson (1980), the \textit{weighted conditional least
squares estimators} (WCLSEs) of the offspring and immigration means of GWI
processes were proposed by Wei and Winnicki (1990), who proved some
self-normalized central limit theorems for the estimators. The limiting
distributions in Wei and Winnicki (1990) are also Gaussian except in the
critical case. The reader can refer to de la Pe\~na et al.\ (2009) for
recent developments on self-normalized limit theorems and their
statistical applications. The estimation problems of the CIR-model defined
by (\ref{1.1}) were studied by Overbeck and Ryd\'{e}n (1997). They
proposed some CLSEs and WCLSEs and proved a Gaussian central limit theorem
for them; see also Overbeck (1998).

In this work, we give some estimation of the drift coefficients $(b,a)$
of the SCIR-model using low frequency observations at equidistant time
points $\{k\Delta: k=0,1,\cdots, n\}$ of a single realization $\{X_t:
t\ge 0\}$. For simplicity, we take $\Delta=1$, but all the results
presented below can be modified to the general case. We shall also
consider the parameters
 \beqlb\label{1.4}
\gamma=e^{-b}, \qquad \rho = ab^{-1}(1-\gamma).
 \eeqlb
Following Klimko and Nelson (1978) and Overbeck and Ryd\'{e}n (1997), we
first define the CLSEs of the parameters. The basic ideas are explained
as follows. By applying It\^{o}'s formula to (\ref{1.2}), for any $t\ge
r\ge 0$ we have
 \beqlb\label{1.5}
X_t = e^{-b(t-r)}X_r + a\int_r^t e^{-b(t-s)}ds + \sigma\int_r^t
e^{-b(t-s)}X_{s-}^{1/\alpha}dZ_s.
 \eeqlb
From (\ref{1.5}) we obtain the stochastic regressive equation
 \beqlb\label{1.6}
X_k = \rho+\gamma X_{k-1}+\varepsilon_k,
 \eeqlb
where
 \beqlb\label{1.7}
\varepsilon_k = \sigma \int_{k-1}^ke^{-b(k-s)}X_{s-}^{1/\alpha} dZ_s.
 \eeqlb
One can see that $\{\varepsilon_k: k\ge 0\}$ is a sequence of martingale
differences with respect to $\{\mcr{F}_k: k\ge 0\}$. The CLSEs of
$(\gamma,\rho)$ and $(b,a)$ can be given by minimizing the sum of squares
 \beqlb\label{1.8}
\sum_{k=1}^n \varepsilon_k^2
 =
\sum_{k=1}^n (X_k-\gamma X_{k-1}-\rho)^2.
 \eeqlb
In particular, those estimators of $(b,a)$ are given by
 \beqlb\label{1.9}
\hat{b}_n=-\log \frac {\sum_{k=1}^nX_{k-1}\sum_{k=1}^nX_k -
n\sum_{k=1}^nX_{k-1}X_k} {\big(\sum_{k=1}^nX_{k-1}\big)^2 -
n\sum_{k=1}^nX_{k-1}^2}
 \eeqlb
and
 \beqlb\label{1.10}
\hat{a}_n=\frac{\sum_{k=1}^nX_{k}-e^{-\hat{b}_n}\sum_{k=1}^n
X_{k-1}}{n(1-e^{-\hat{b}_n})}\hat{b}_n.
 \eeqlb
Following Wei and Winnicki (1990), we also consider the WCLSEs of
$(\gamma,\rho)$ and $(b,a)$ by minimizing the weighted sum
 \beqlb\label{1.11}
\sum_{k=1}^n \frac{\varepsilon_k^2}{X_{k-1} + 1}
 =
\sum_{k=1}^n \frac{[X_k - \gamma(X_{k-1}+1) - (\rho-\gamma)]^2} {X_{k-1}
+ 1}.
 \eeqlb
The reason of considering the above quantity is it does not fluctuate too
much even when the values of the samples $X_k, k=0,1,\cdots, n$ are
large. The resulting WCLSEs of $(b,a)$ are given by
 \beqlb\label{1.12}
\check{b}_n=-\log\frac{\sum_{k=1}^nX_k\sum_{k=1}^n\frac{1}{X_{k-1}+1}-n
\sum_{k=1}^n\frac{X_k}{X_{k-1}+1}}
{\sum_{k=1}^n(X_{k-1}+1)\sum_{k=1}^n\frac{1}{X_{k-1}+1}-n^2}
 \eeqlb
and
 \beqlb\label{1.13}
\check{a}_n=\frac{\sum_{k=1}^nX_{k}-e^{-\check{b}_n}\sum_{k=1}^n
X_{k-1}}{n(1-e^{-\check{b}_n})}\check{b}_n.
 \eeqlb

The main purpose of this paper is to study the asymptotic properties of
the CLSEs and the WCLSEs given above. We show that the estimators are
consistent and obey some central limit theorems. In particular, for
$1<\alpha\le 2$ we prove that $n^{(\alpha-1)/\alpha} (\check{b}_n-b,
\check{a}_n-a)$ converges to an $\alpha$-stable random vector as $n\to
\infty$. For $1<\alpha<(1+\sqrt{5})/2$, we show that
$n^{(\alpha-1)/\alpha^2} (\hat{b}_n-b, \hat{a}_n-a)$ converges to a
nontrivial limit as $n\to \infty$. A combination of this with the result
of Overbeck and Ryd\'{e}n (1997) for $\alpha=2$ only gives a partial
characterization of the asymptotic behavior of the CLSEs. The
characterization of the exact asymptotics of the CLSEs for
$(1+\sqrt{5})/2\le \alpha< 2$ is left as an open problem.

The proofs of our limit theorems are very different from and much harder
than the Gaussian case. The key of the approach is to establish the
convergence of some point processes and partial sums associated with a
stationary realization of the SCIR-model. The techniques in this subject
have been developed extensively by Basrak and Segers (2009), Davis and
Hsing (1995), Davis and Mikosch (1998) among others. We also make use of
the results of Hult and Lindskog (2007) on the extremal behavior of
L\'{e}vy stochastic integrals. The proofs depend heavily on the
construction and characterization of CBI-processes in terms of stochastic
equations of given in Dawson and Li (2006, 2012), Fu and Li (2010) and Li
and Ma (2008).

We finally propose an estimator of the volatility coefficient $\sigma$
based on high frequency observations at times $\{0,1/n,\cdots,
(n-1)/n,1\}$. Suppose that the parameter $\alpha$ is known. Given
constants $p\in (0,\alpha)$ and $\delta\in (0, \min\{1-1/\alpha,
1/\alpha^2\})$, let
 \beqlb\label{1.14}
\hat{\sigma}_n=\frac{1}{n^{1/p-1/\alpha}\mbf{E}^{1/p}[|Z_1|^p]}
\Big(\sum_{k=1}^n\Big|\frac{X_{k/n}-X_{(k-1)/n}} {X_{(k-1)/n}^{1/\alpha} +
n^{-\delta}}\Big|^p\Big)^{1/p}.
 \eeqlb
We prove that $\hat{\sigma}_n$ is a weakly consistent estimator for
$\sigma$.

The paper is organized as follows. In Section 2, we prove the exponential
ergodicity of some subcritical CBI-processes, which implies the strong
mixing property of the SCIR-model. Section 3 is devoted to the regular
variation properties of some random sequences defined from the model. The
limit theorems of random point processes and partial sums are established
in Section 4. Based on those theorems, the asymptotic properties of the
estimators are proved in Section 5.

\medskip

\textbf{Notation.}~ Let $\mbb{N} = \{0, 1, 2, \cdots\}$ and $\mbb{Z} =
\{0, \pm 1, \pm 2, \cdots\}$. Let $\mbb{R} = (-\infty,\infty)$,
$\bar{\mbb{R}} = [-\infty,\infty]$ and $\bar{\mbb{R}}_0^d =
\bar{\mbb{R}}^d\setminus \{\mbf{0}\}$, where $\mbf{0} = (0,0,\cdots,0)$.
Let $C_0^+(\bar{\mbb{R}}_0^d)$ be the collection of positive continuous
functions on $\bar{\mbb{R}}_0^2$ with compact support. Let
$M(\bar{\mbb{R}}_0^d)$ be the class of Radon point measures on
$\bar{\mbb{R}}_0^d$ furnished with the topology of vague convergence. We
use $C$ with or without subscripts to denote positive constants whose
values are not important.


\section{CBI-processes and ergodicity}

\setcounter{equation}{0}

In this section, we prove some simple properties of CBI-processes. In
particular, we prove a subcritical CBI-process is exponentially ergodic
and strongly mixing. The results are useful in the study of the
asymptotics of the estimators. We start with an important special case of
those processes. Let $\sigma\ge 0$ and $b$ be constants and $(u\wedge
u^2)m(du)$ a finite measure on $(0,\infty)$. For $z\ge 0$ set
 \beqnn
\phi(z) = bz + \frac{1}{2}\sigma^2 z^2 + \int_0^\infty(e^{-zu}-1+zu)m(du).
 \eeqnn
A Markov process with state space $\mbb{R}_+ := [0,\infty)$ is called a
\textit{continuous-state branching process} (CB-process) with
\textit{branching mechanism} $\phi$ if it has transition semigroup
$(Q_t)_{t\ge 0}$ given by
 \beqlb\label{2.1}
\int_0^\infty e^{-\lambda y}Q_t(x,dy)
 =
e^{-x v_t(\lambda)},
 \eeqlb
where $t\mapsto v_t(\lambda)$ is the unique positive solution of
 \beqlb\label{2.2}
\frac{\partial}{\partial t}v_t(\lambda) = -\phi(v_t(\lambda)), \qquad
v_0(\lambda)=\lambda.
 \eeqlb
The CB-process is called \textit{critical}, \textit{subcritical} or
\textit{supercritical} as $b=0$, $b>0$ or $b<0$, respectively. From
(\ref{2.2}) we obtain the following semigroup property:
 \beqnn
v_{r+t}(\lambda) = v_r(v_t(\lambda)), \qquad r,t,\lambda\ge 0.
 \eeqnn
Taking the derivatives of both sides of (\ref{2.2}) one can see $u_t :=
(d/d\lambda) v_t(0)$ solves the equation $(d/dt) u_t = -bu_t$, and so $u_t
= e^{-bt}$ for $t\ge 0$. Then differentiating both sides of (\ref{2.1})
gives
 \beqlb\label{2.3}
\int_0^\infty y Q_t(x,dy)
 =
xe^{-bt}, \qquad t,x\ge 0.
 \eeqlb
By Jensen's inequality, we have $v_t(\lambda)\le \lambda e^{-bt}$ for
$t,\lambda\ge 0$.

It is easy to see that $(Q_t)_{t\ge 0}$ is a Feller transition semigroup,
so it has a Hunt realization. Let $X = (\Omega, \mcr{G}, \mcr{G}_t, X_t,
\mbf{Q}_x)$ be a Hunt realization of the CB-process. The hitting time
$\tau_0 = \inf\{t\ge 0: X_t=0\}$ is called the \textit{extinction time}
of $X$. It follows from Theorem~3.5 of Li (2011) that for $t\ge 0$ the
limit $\bar{v}_t = \uparrow \lim_{\lambda\to \infty} v_t(\lambda)$ exists
in $(0,\infty]$, and
 \beqlb\label{2.4}
\mbf{Q}_x(\tau_0\le t) = \mbf{Q}_x(X_t=0) = \exp\{-x\bar{v}_t\}.
 \eeqlb
By Theorem~3.8 of Li (2011), we have $\bar{v}_t< \infty$ for all $t>0$ if
and only if the following condition holds:

\bcondition\label{t2.1} There is some constant $\theta> 0$ such that
$\phi(z)>0$ for $z> \theta$ and
 \beqnn
\int_{\theta}^\infty \phi(z)^{-1}dz< \infty.
 \eeqnn
\econdition

\bproposition\label{t2.2} (Li, 2011, p.61) If Condition~\ref{t2.1} holds,
then $\bar{v}_t = l_t(0,\infty)$ is the minimal solution of
 \beqlb\label{2.5}
\frac{d}{dt}\bar{v}_t = -\phi(\bar{v}_t), \qquad \bar{v}_0 = \infty.
 \eeqlb
\eproposition

Let $t\mapsto v_t(\lambda)$ be defined by (\ref{2.2}). A Markov process
with state space $\mbb{R}_+$ is called a \textit{CBI-process} with
\textit{branching mechanism} $\phi$ and \textit{immigration rate} $a\ge
0$ if it has transition semigroup $(P_t)_{t\ge 0}$ given by
 \beqlb\label{2.6}
\int_0^\infty e^{-\lambda y}P_t(x,dy)
 =
\exp\Big\{-xv_t(\lambda) - a\int_0^tv_s(\lambda)ds\Big\}.
 \eeqlb
By differentiating both sides of (\ref{2.6}) we obtain
 \beqlb\label{2.7}
\ar\ar~ \int_0^\infty y P_t(x,dy)
 =
xe^{-bt} + a\int_0^t e^{-bs} ds
 =
xe^{-bt} + ab^{-1}(1-e^{-bt}),
 \eeqlb
where $b^{-1}(1-e^{-bt}) = t$ when $b=0$ by convention.

A realization of the CBI-process can be constructed as the strong
solution to a stochastic integral equation. Let $W(ds,du)$ be a
time-space Gaussian white noise on $(0,\infty)^2$ with intensity $dsdu$
and $N_1(ds,dz,du)$ a Poisson random measure on $(0,\infty)^3$ with
intensity $dsm(dz)du$. Let $\tilde{N}_1(ds,dz,du) = N_1(ds,dz,du) -
dsm(dz)du$ denote the compensated measure. Then for each $x\ge 0$ there
is a pathwise unique positive strong solution to the following stochastic
equation:
 \beqlb\label{2.8}
Y_t(x)
 \ar=\ar
x + \int_0^t(a-bY_s(x)) ds + \sigma\int_0^t\int_0^{Y_{s-}(x)} W(ds,du) \cr
 \ar\ar\qquad
+ \int_0^t\int_0^\infty\int_0^{Y_{s-}(x)} z \tilde{N}_1(ds,dz,du).
 \eeqlb
The solution $\{Y_t(x),t\ge 0\}$ is a CBI-process with branching
mechanism $\phi$ and immigration rate $a$. See Theorem~3.1 of Dawson and
Li (2012) or Theorem~2.1 of Li and Ma (2008). A slightly different
formulation of the process was given in Dawson and Li (2006).

\bproposition\label{t2.3} Suppose that Condition~\ref{t2.1} holds. For
$x,y\ge 0$ let $T_{x,y} :=\inf\{t\ge 0: Y_t(x) - Y_t(y) = 0\}$. Then we
have $\mbf{P}\{T_{x,y}< \infty\} = 1$ and
 \beqlb\label{2.9}
\mbf{P}\{T_{x,y}\le t\}=\exp\{-|x-y|\bar{v}_t\}, \qquad t\ge 0.
 \eeqlb
Moreover, we have $Y_t(x) = Y_t(y)$ for all $t\ge T_{x,y}$. \eproposition

\proof It suffices to consider the case of $y\ge x\ge 0$. By Theorem~3.2
of Dawson and Li (2012), we have $\mbf{P}\{Y_t(x)\ge Y_t(y)\ge 0$ for all
$t\ge 0\}=1$ and $\{Y_t(x)-Y_t(y): t\ge 0\}$ is a CB-process with
branching mechanism $\phi$; see also Remark~2.1~(iv) of Li and Ma (2008).
Then (\ref{2.9}) follows from (\ref{2.4}). The pathwise uniqueness of
(\ref{2.8}) implies that $Y_t(x)=Y_t(y)$ for all $t\ge T_{x,y}$. By
Corollary~3.9 of Li (2011) we have $\mbf{P}\{T_{x,y}< \infty\} = 1$. \qed

The above proposition provides a successful coupling of the CBI-processes.
This has many important implications. We refer the reader to Chen (2004)
for systematical study of coupling methods and their applications in the
theory of Markov processes. In particular, we shall use the above coupling
to prove the strong Feller property and exponential ergodicity of the
CBI-process following Chen (2004, p.37). Write $f\in \b\mcr{B}(\mbb{R}_+)$
if $f$ is a bounded measurable function on $\mbb{R}_+$.

\ttheorem\label{t2.4} Under Condition~\ref{t2.1}, the transition semigroup
$(P_t)_{t\ge 0}$ given by (\ref{2.6}) has the strong Feller property.
Moreover, for any $t>0$, $x,y\ge 0$ and $f\in \b\mcr{B}(\mbb{R}_+)$ we
have
 \beqlb\label{2.10}
\big|P_tf(x)-P_tf(y)\big|\le 2\|f\|_\infty (1-e^{-\bar{v}_t|x-y|}),
 \eeqlb
where $\|f\|_\infty = \sup_x|f(x)|$ denotes the supremum norm. \etheorem

\proof It suffices to prove (\ref{2.10}). By Proposition~\ref{t2.3} we
have
 \beqnn
|P_tf(x)-P_tf(y)| = |\mbf{E}[f(Y_t(x))-f({Y}_t(y))]|
 \le
2\|f\|_\infty \mbf{P}(T_{x,y}>t),
 \eeqnn
which together with (\ref{2.9}) implies (\ref{2.10}). \qed

\ttheorem\label{t2.5} Suppose that $b>0$. Then the transition semigroup
$(P_t)_{t\ge 0}$ has a unique stationary distribution $\mu$, which is
given by
 \beqlb\label{2.11}
\ar\ar~ L_{\mu}(\lambda)
 =
\int_0^\infty e^{-\lambda x} \mu(dx)
 =
\exp\Big\{-a\int_0^\lambda z\phi(z)^{-1}dz\Big\},
 ~~
\lambda\ge 0.
 \eeqlb
Moreover, we have
 \beqlb\label{2.12}
\int_0^\infty x\mu(dx)
 =
\frac{d}{d\lambda}L_\mu(\lambda)\Big|_{\lambda=0+} = \frac{a}{b}.
 \eeqlb
\etheorem

\proof By Theorem~3.20 of Li (2011) and its proof given there, for any
$x\ge 0$ we have $\mu = \lim_{t\to \infty} P_t(x,\cdot)$ by the weak
convergence. Then $(P_t)_{t\ge 0}$ has the unique stationary distribution
$\mu$. By differentiating (\ref{2.11}) we obtain (\ref{2.12}). \qed

\ttheorem\label{t2.6} Suppose that $b>0$ and Condition~\ref{t2.1} is
satisfied. Then the transition semigroup $(P_t)_{t\ge 0}$ is exponentially
ergodic. More precisely, for any $x\ge 0$ and $t\ge 1$ we have
 \beqlb\label{2.13}
\ar\ar~ \|P_t(x,\cdot)-\mu(\cdot)\|_{\rm var}
 \le
2(1-\exp\{-\bar{v}_1xe^{-b(t-1)}\}) + 2\bar{v}_1ab^{-1}e^{-b(t-1)},
 \eeqlb
where $\mu$ is given by (\ref{2.11}) and $\|\cdot\|_{\rm var}$ denotes
the total variation norm. \etheorem

\proof We only need to prove (\ref{2.13}). In view of (\ref{2.11}), we
have
 \beqlb\label{2.14}
0\le 1-L_\mu(\lambda)\le a\int_0^\lambda z\phi(z)^{-1}dz
 \le
ab^{-1}\lambda.
 \eeqlb
Then, for any $f\in \b\mcr{B}(\mbb{R}_+)$ with $\|f\|\le 1$, we can use
Theorem~\ref{t2.4} and (\ref{2.14}) to see
 \beqnn
|P_tf(x)-\mu(f)|
 \ar\le\ar
\int_0^\infty|P_tf(x)-P_tf(y)|\mu(dy)
 \le
2\int_0^\infty(1-e^{-\bar{v}_t|x-y|})\mu(dy) \cr
 \ar\le\ar
2\int_0^x (1-e^{-\bar{v}_tx}) \mu(dy) + 2\int_x^\infty
(1-e^{-\bar{v}_ty}) \mu(dy) \ccr
 \ar\le\ar
2(1-e^{-\bar{v}_tx}) + 2[1-L_\mu(\bar{v}_t)]
 \le
2(1-e^{-\bar{v}_tx}) + 2ab^{-1}\bar{v}_t.
 \eeqnn
By Proposition~\ref{t2.2}, for $t\ge 1$ we have $\bar{v}_t =
v_{t-1}(\bar{v}_1)$ and so $\bar{v}_t\le e^{-b(t-1)}\bar{v}_1$. Then we
obtain (\ref{2.13}). \qed

Under the conditions of Theorem~\ref{t2.6}, for any finite set $\{t_1<
t_2< \cdots< t_n\}\subset \mbb{R}$ we can define the probability measure
$\mu_{t_1,t_2,\cdots,t_n}$ on $\mbb{R}_+^n$ by
 \beqlb\label{2.15}
\ar\ar \mu_{t_1,t_2,\cdots,t_n}(dx_1,dx_2,\cdots,dx_n) \ccr
 \ar\ar\quad
= \mu(dx_1)P_{t_2-t_1}(x_1,dx_2)\cdots P_{t_n-t_{n-1}}(x_{n-1},dx_n).
 \eeqlb
It is easy to see that $\{\mu_{t_1,t_2,\cdots,t_n}: t_1< t_2< \cdots<
t_n\in \mbb{R}\}$ is a consistent family. By Kolmogorov's theorem, there
is a stochastic process $\{Y_t: t\in \mbb{R}\}$ with finite-dimensional
distributions given by (\ref{2.15}). This process is a (strictly)
stationary Markov process with one-dimensional marginal distribution $\mu$
and transition semigroup $(P_t)_{t\ge 0}$. Since $(P_t)_{t\ge 0}$ is a
Feller semigroup, the process $\{Y_t: t\in \mbb{R}\}$ has a c\`{a}dl\`{a}g
modification.

\ttheorem\label{t2.7} Let $\{Y_t: t\in \mbb{R}\}$ be a Markov process with
finite-dimensional distributions given by (\ref{2.15}). Then it is
strongly mixing with geometric rate, that is, as $t\to \infty$,
 \beqnn
\pi_t := \sup_{A\in\sigma\{Y_s,s\le0\}}\sup_{B\in\sigma\{Y_s,s>t\}}
\big|\mbf{P}(A\cap B)-\mbf{P}(A)\mbf{P}(B)\big|
 \eeqnn
decays to zero exponentially. \etheorem

\proof It follows from (\ref{2.13}) and (\ref{2.14}) that, for $t\ge 1$,
 \beqnn
\int_0^\infty\|P_t(x,\cdot)-\mu(\cdot)\|_{\rm var}\mu(dx)
 \ar\le\ar
2\int_0^\infty (1-e^{-\bar{v}_1xe^{-b(t-1)}})\mu(dx) \cr
 \ar\ar\qquad\qquad\qquad
+\, 2\bar{v}_1ab^{-1} e^{-b(t-1)} \ccr
 \ar\le\ar
2[1-L_\mu(\bar{v}_1e^{-b(t-1)})] + 2\bar{v}_1ab^{-1} e^{-b(t-1)} \ccr
 \ar\le\ar
4\bar{v}_1ab^{-1} e^{-b(t-1)}.
 \eeqnn
Then $\{Y_t: t\in \mbb{R}\}$ is strongly mixing with geometric rate; see,
e.g., Mikosch and Straumann (2006, p.516) or Bradley (2005, p.112). \qed

\blemma\label{t2.8} Suppose that $\{Z_t\}$ is an $\alpha$-stable L\'evy
process with $1<\alpha<2$ and $\{y(t)\}$ is a predictable process
satisfying, a.s.,
 \beqnn
\int_0^T |y(t)|^\alpha dt< \infty, \qquad T\ge 0.
 \eeqnn
Then for any $0<r<\alpha$, there exists a constant $C=C(r,\alpha)\ge 0$
such that
 \beqnn
\mbf{E}\Big[\sup_{t\le T}\Big|\int_0^ty(s)dZ_s\Big|^r\Big]
 \le
C\mbf{E}\Big[\Big(\int_0^T |y(t)|^\alpha dt\Big)^{r/\alpha}\Big].
 \eeqnn
\elemma

\proof When $\{y(t)\}$ is a positive process, the result follows from a
result proved in Long and Qian (2011). Then, in the general case, we have
 \beqnn
\mbf{E}\Big[\sup_{t\le T}\Big|\int_0^ty(s)dZ_s\Big|^r\Big]
 \ar\le\ar
C_1\mbf{E}\Big[\sup_{t\le T}\Big|\int_0^ty_+(s)dZ_s\Big|^r + \sup_{t\le
T}\Big|\int_0^ty_-(s)dZ_s\Big|^r\Big] \cr
 \ar\le\ar
C_2\mbf{E}\Big[\Big(\int_0^T y_+(t)^{\alpha} dt\Big)^{r/\alpha} +
\Big(\int_0^T y_-(t)^{\alpha} dt\Big)^{r/\alpha}\Big] \cr
 \ar\le\ar
2C_2\mbf{E}\Big[\Big(\int_0^T |y(t)|^{\alpha} dt\Big)^{r/\alpha}\Big],
 \eeqnn
where $y_+$ and $y_-$ denote respectively the positive and negative parts
of $y$. \qed

Now let us consider a filtered probability space $(\Omega,\mcr{F},
\mcr{F}_t,\mbf{P})$ satisfying the usual hypotheses. Let $\{Z_t: t\ge
0\}$ be a spectrally positive $\alpha$-stable L\'evy process. For
$\alpha=2$ we understand the process as a standard Brownian motion; and
for $1<\alpha<2$ we assume it is a stable process with L\'evy measure
$\nu_\alpha(dz)$ given by (\ref{1.3}). By Theorem~6.2 of Fu and Li
(2010), for any initial value $X_0$, which is a positive
$\mcr{F}_0$-measurable random variable, there is a unique positive strong
solution $\{X_t: t\ge 0\}$ to (\ref{1.2}). The existence and uniqueness
of this solution also follows from Corollary~6.3 of Fu and Li (2010) by a
time change. Let $f$ be a bounded continuous function on $\mbb{R}$ with
bounded continuous derivatives up to the second order. For $\alpha=2$, we
can use It\^o's formula to see that
 \beqlb\label{2.16}
f(X_t) = f(X_r) + \int_r^t Lf(X_s) ds + M_t(f), \qquad t\ge r,
 \eeqlb
where $\{M_t(f): t\ge r\}$ is a martingale with respect to the filtration
$(\mcr{F}_t)_{t\ge r}$ and
 \beqnn
Lf(x) = (a-bx)f'(x) + \frac{\sigma^2}{2} xf''(x), \qquad x\ge 0.
 \eeqnn
When $1<\alpha<2$, by the L\'{e}vy-It\^{o} representation of $\{Z_t\}$,
we can rewrite ($\ref{1.2}$) into the integral form:
 \beqlb\label{2.17}
\ar\ar~ X_t = X_r + \int_r^t(a-bX_s)ds + \sigma \int_r^t\int_0^\infty
X_{s-}^{1/\alpha}z \tilde{N}(ds,dz), \quad t\ge r,
 \eeqlb
where $\tilde{N}(ds,dz)$ is a compensated Poisson random measure on
$(0,\infty)^2$ with intensity $ds\nu_\alpha(dz)$. By It\^o's formula one
can see that (\ref{2.16}) still holds for $1<\alpha<2$ with the operator
$L$ defined by
 \beqnn
Lf(x) = (a-bx)f'(x) + \frac{\sigma^\alpha}{\alpha\Gamma(-\alpha)}
\int_0^\infty [f(x+y) - f(x) - yf'(x)] \frac{dy}{y^{\alpha+1}}.
 \eeqnn
By Theorem~9.30 of Li (2011), for any $1<\alpha\le 2$ we can identify the
SCIR-model as a subcritical CBI-process with immigration rate $a$ and
branching mechanism
 \beqnn
\phi(\lambda) = b\lambda + \frac{\sigma^\alpha}{\alpha} \lambda^\alpha,
\qquad \lambda\ge 0.
 \eeqnn
It follows from Theorem~\ref{t2.5} that the SCIR-model has the unique
stationary distribution $\mu$ with Laplace transform given by
 \beqlb\label{2.19}
\ar\ar~ L_\mu(\lambda) = \int_0^\infty e^{-\lambda x} \mu(dx)
 =
\exp\Big\{-\int_0^\lambda \frac{\alpha adz}{\alpha b+ \sigma^\alpha
z^{\alpha-1}}\Big\}, \quad \lambda\ge 0.
 \eeqlb
Let $\mbf{P}_x$ denote the law of the SCIR-model $\{X_t: t\ge 0\}$ defined
by (\ref{1.2}) with $X_0=x\ge 0$ and let $\mbf{E}_x$ denote the
corresponding expectation.

\bproposition\label{t2.9} Suppose that $0<\alpha<2$. Then for any
$0<\beta<\alpha$, there is a constant $C\ge 0$ so that, for $t,T\ge 0$,
 \beqnn
\mbf{E}_x\Big(\Big|\int_0^t e^{-b(t-s)} X_{s-}^{1/\alpha}
dZ_s\Big|^\beta\Big)
 \le
C(1+x^{\beta/\alpha}e^{-\beta bt/\alpha})
 \eeqnn
 and
 \beqnn
\mbf{E}_x\Big(\sup_{0\le t\le T}\Big|\int_0^t e^{-b(t-s)}
X_{s-}^{1/\alpha} dZ_s\Big|^\beta\Big)
 \le
C(x^{\beta/\alpha}e^{\beta b(1-1/\alpha)T} + e^{\beta bT}).
 \eeqnn
\eproposition

\proof Using Lemma~\ref{t2.8}, we have
 \beqnn
\mbf{E}_x\Big(\Big|\int_0^t e^{-b(T-s)} X_{s-}^{1/\alpha}
dZ_s\Big|^\beta\Big)
 \le
C_1e^{-\beta bt}\mbf{E}_x\Big[\Big(\int_0^t e^{\alpha bs} X_s
ds\Big)^{\beta/\alpha}\Big]
 \eeqnn
and
 \beqnn
\mbf{E}_x\Big(\sup_{0\le t\le T}\Big|\int_0^t e^{-b(t-s)}
X_{s-}^{1/\alpha} dZ_s\Big|^\beta\Big)
 \le
C_1\mbf{E}_x\Big[\Big(\int_0^T e^{\alpha bs} X_s
ds\Big)^{\beta/\alpha}\Big].
 \eeqnn
By H\"older's inequality and (\ref{2.7}) it is easy to see
 \beqnn
\mbf{E}_x\Big[\Big(\int_0^t e^{\alpha bs} X_s ds\Big)^{\beta/\alpha}\Big]
 \ar\le\ar
\Big[\int_0^t\mbf{E}_x(e^{\alpha bs}X_s) ds\Big]^{\beta/\alpha} \cr
 \ar\le\ar
\Big[\int_0^t e^{\alpha bs}(xe^{-bs} + ab^{-1})ds\Big]^{\beta/\alpha} \ccr
 \ar\le\ar
C_2(x^{\beta/\alpha}e^{\beta b(1-1/\alpha)t} + e^{\beta bt}).
 \eeqnn
Then we have the desired inequalities. \qed

\bproposition\label{t2.10} Suppose that $0<\alpha<2$. Then for any
$0<\beta<\alpha$, there is a constant $C\ge 0$ and a locally bounded
function $T\mapsto C(T)\ge 0$ so that, for $t,T\ge 0$,
 \beqnn
\mbf{E}_x(X_t^\beta)\le C(1+x^\beta e^{-\beta bt/\alpha})
 \eeqnn
and
 \beqnn
\mbf{E}_x\Big(\sup_{0\le t\le T}X_t^\beta\Big)\le C(T)(1+x^\beta).
 \eeqnn
\eproposition

\proof Using (\ref{1.5}) with $r=0$ and an elementary inequality, we have
 \beqnn
\mbf{E}_x(X_t^\beta)
 \le
C_1\mbf{E}_x\Big[x^\beta e^{-\beta bt} + a^\beta b^{-\beta} +
\sigma^\beta\mbf{E}_x \Big(\Big|\int_0^t e^{-b(t-s)}X_{s-}^{1/\alpha}
dZ_s\Big|^\beta\Big)\Big]
 \eeqnn
and
 \beqnn
\mbf{E}_x\Big(\sup_{0\le t\le T}X_t^\beta\Big)
 \ar\le\ar
C_1\Big[x^\beta + a^\beta b^{-\beta} +
\sigma^\beta\mbf{E}_x\Big(\sup_{0\le t\le T} \Big|\int_0^t
e^{-b(t-s)}X_{s-}^{1/\alpha} dZ_s\Big|^\beta\Big)\Big].
 \eeqnn
Then the results follow by Proposition~\ref{t2.9}. \qed


\section{Regular variations}

\setcounter{equation}{0}

In this section, we study the regular variation property of some random
sequences associated with the SCIR-model. We shall make use of the
stochastic equations (\ref{1.2}) and (\ref{1.5}) as well as the results of
Hult and Lindskog (2007) on the extremal behavior of L\'{e}vy stochastic
integrals. The reader may also refer to Samorodnitsky and Grigoriu (2003)
for some results on the tail behavior of solutions to certain stochastic
differential equations driven by L\'evy processes. For the convenience of
the reader, we first recall some concepts and properties of regularly
varying sequences. Let ``$|\cdot|$'' be any norm on $\mbb{R}^d$.

\bdefinition\label{t3.1} A $d$-dimensional random vector $\mbf{X}$ is
said to be \textit{regularly varying} if there exists a Radon measure
$\eta$ on $\mbb{R}^d$, finite on sets of the form $\{x\in \mbb{R}^d:
|x|\ge r\}$, and a sequence $\{a_n\}$ satisfying $a_n\to \infty$ such
that, as $n\to \infty$,
 \beqlb\label{3.1}
n\mbf{P}(a_n^{-1}\mbf{X}\in \cdot)\overset{\rm v}{\longrightarrow}
\eta(\cdot).
 \eeqlb
The above sequential form of the condition is the same as saying there
exists a (necessarily regularly varying) function $t\mapsto g(t)$ such
that, as $n\to \infty$,
 \beqnn
g(t)\mbf{P}(t^{-1}\mbf{X}\in \cdot)\overset{\rm v}{\longrightarrow}
\eta(\cdot).
 \eeqnn
It is known that the condition implies the existence of a constant
$\alpha>0$ such that $\eta(rA) = r^{-\alpha}\eta(A)$ for all $r>0$ and
all $A\in \mcr{B}(\mbb{R}^d)$ bounded away from $\bf 0$, where $rA=\{rx:
x\in A\}$. In this case, we say $\mbf{X}$ is regularly varying
\textit{with index $\alpha>0$}. \edefinition

\bproposition\label{t3.2} Let $\xi$ be a positive regularly varying
random variable with index $\alpha>0$. Then we have:
 \benumerate

\itm[\rm(i)] If $\alpha>1$, then as $x\to \infty$,
 \beqnn
\mbf{E}(\xi 1_{\{\xi>x\}})
 \sim
\frac{\alpha}{\alpha-1} x\mbf{P}(\xi>x);
 \eeqnn

\itm[\rm(ii)] If $0<\alpha<1$, then as $x\to \infty$,
 \beqnn
\mbf{E}(\xi 1_{\{\xi<x\}})
 \sim
\frac{\alpha}{1-\alpha} x\mbf{P}(\xi>x).
 \eeqnn

 \eenumerate
\eproposition

\proof The results are immediate consequences of Karamata's theorem; see
Resnick (2007, p.25). Let $G(x) = \mbf{P}(\xi>x)$. By integration by
parts,
 \beqnn
\mbf{E}(\xi 1_{\{\xi>x\}})
 =
-\int_x^\infty y dG(y)
 =
x G(x) + \int_x^\infty G(y)dy.
 \eeqnn
The second term on the right hand side is equivalent to $(\alpha-1)^{-1}
xG(x)$. Then (i) follows. The proof of (ii) is similar; see also Resnick
(2007, p.36). \qed

There are several equivalent forms of the regular variation property;
see, e.g., Resnick (1986, p.69). One of them is given in the next
theorem. Let $\mbb{S}^{d-1} = \{x\in \mbb{R}^d: |x|=1\}$ be the unit
sphere.

\ttheorem\label{t3.3} A $d$-dimensional random vector $\mbf{X}$ is
regularly varying with index $\alpha>0$ if and only if there exists a
probability measure $\nu$ on $\mbb{S}^{d-1}$ such that for every $r>0$,
as $t\to \infty$,
 \beqnn
\frac{\mbf{P}(|\mbf{X}|>rt, |\mbf{X}|^{-1}\mbf{X}\in \cdot)}
{\mbf{P}(|\mbf{X}|> t)}
 \overset{\rm w}{\longrightarrow}
r^{-\alpha}\nu(\cdot).
 \eeqnn
\etheorem

\bdefinition\label{t3.4} A sequence of random variables $\{\mbf{X}_k:
k\in\mbb{Z}\}$ in $\mbb{R}^d$ is called \textit{jointly regularly
varying} if all the vectors of the form $(\mbf{X}_1,\cdots,\mbf{X}_l)$
are regularly varying. \edefinition

In the sequel, we shall also use the norm $\|\mbf{x}\| := \max_i|x_i|$
for $\mbf{x}=(x_1,\cdots,x_d)\in \mbb{R}^d$. By Corollary~3.2 in Basrak
and Segers (2009) we have the following equivalent characterization of
the jointly regular variation property:

\ttheorem\label{t3.5} Let $\{\mbf{X}_k: k\in\mbb{Z}\}$ be a stationary
sequence in $\mbb{R}^d$ and assume $x\mapsto \mbf{P}(\|\mbf{X}_0\|>x)$ is
regularly varying with index $-\alpha$ for some $\alpha>0$. Then
$\{\mbf{X}_k: k\in\mbb{Z}\}$ is jointly regularly varying with index
$\alpha$ if and only if there is a sequence $\{\mbf{\Theta}_k: k\in
\mbb{N}\}$ such that for every $k\in\mbb{N}$, as $x\to \infty$,
 \beqnn
\mbf{P}\big(\|\mbf{X}_0\|^{-1}(\mbf{X}_0,\cdots,\mbf{X}_k)\in \cdot\big|
\|\mbf{X}_0\|>x\big)
 \overset{\rm w}{\longrightarrow}
\mbf{P}\big((\mbf{\Theta}_0,\cdots, \mbf{\Theta}_k)\in\cdot\big).
 \eeqnn
\etheorem

The concept of regular variations can also be defined for continuous time
stochastic processes. Let $T\ge 0$ and let $\mbb{D}^d[0,T] :=
\mbb{D}([0,T],\mbb{R}^d)$ be the space of all $\mbb{R}^d$-valued
c\`{a}dl\`{a}g functions on $[0,T]$ equipped with Skorokhod topology. Let
 \beqnn
\mbb{S}^d[0,T] = \Big\{{\bf{y}}\in \mbb{D}^d[0,T]: \sup_{0\le t\le
T}\|\mbf{y_t}\| = 1\Big\}.
 \eeqnn

\bdefinition\label{t3.6} A stochastic process $\mbf{Y} = \{\mbf{Y}_t:
0\le t\le T\}$ with sample path in $\mbb{D}^d[0,T]$ is said to be
\textit{regularly varying} if there exist a measure $Q$ on
$\mbb{D}^d[0,T]$, finite on sets bounded away from $\mbf{0}$, and a
sequence $\{a_n\}$ satisfying $a_n\to \infty$ such that for any set $B\in
\mcr{B}(\mbb{D}^d[0,T])$ bounded away from $\bf 0$ with $Q(\partial
B)=0$, as $n\to \infty$,
 \beqnn
n\mbf{P}(a_n^{-1}\mbf{Y}\in B)\longrightarrow Q(B),
 \eeqnn
The above property implies there is a constant $\alpha>0$ such that
$Q(uB) = u^{-\alpha}Q(B)$ for all $u>0$ and all $B\in
\mcr{B}(\mbb{D}^d[0,T])$ bounded away from $\bf 0$. In this situation, we
say $\mbf{Y}$ is regularly varying \textit{with index $\alpha>0$}
\edefinition

The convergence in the above definition can be formulated for general
boundedly finite measures on $\bar{\mbb{D}}_0^d[0,T] = (0,\infty]\times
\mbb{S}^d[0,T]$. We shall denote the convergence by ``$\overset{\rm
\hat{w}}{\longrightarrow}$''. The reader may refer to Hult and Lindskog
(2005) for more details.

\bremark\label{t3.7} Let $0<\alpha<2$ and let $\{Z_t: t\ge 0\}$ be a
one-dimensional $\alpha$-stable process with L\'evy measure $\nu(dz)$. It
follows from Lemma~2.1 of Hult and Lindskog (2007) that, as $n\to
\infty$,
 \beqnn
n\mbf{P}(n^{-1/\alpha}Z_t\in\cdot)
 \overset{\rm v}{\longrightarrow}
t\nu(\cdot).
 \eeqnn
\eremark

\bremark\label{t3.8} Let $0<\alpha<2$ and $T\ge 1$. Suppose that $\{Z_t:
0\le t\le T\}$ is a one-dimensional L\'evy process such that $X=Z_1$
satisfies (\ref{3.1}) with $\eta(z,\infty) = cz^{-\alpha}$ for some
$c>0$. Let $\{Y_t: 0\le t\le T\}$ be a positive predictable c\`agl\`ad
process satisfying $\sup_{0\le t\le T}Y_t> 0$ a.s.\ and
$\mbf{E}[\sup_{0\le t\le T}Y_t^{\alpha+\delta}]< \infty$ for some
$\delta>0$. By Theorem~3.4 and Example~3.1 in Hult and Lindskog (2007),
for any $z>0$ and $0\le t\le T$ we have, as $n\to \infty$,
 \beqnn
n\mbf{P}\Big(a_n^{-1}\int_0^t Y_{s-} dZ_s> z\Big)
 \to
\eta(z,\infty)\int_0^t \mbf{E}(Y_s^\alpha) ds
 =
cz^{-\alpha}\int_0^t \mbf{E}(Y_s^\alpha) ds.
 \eeqnn
\eremark

\bremark\label{t3.9} Suppose that $\{\mbf{Y}_k\}$ is a stationary sequence
of regularly varying random vectors. Let $\{a_n\}$ be taken such that
$n\mbf{P}(|\mbf{Y}_0|> a_n)\to 1$ as $n\to \infty$. By Lemma~2.3.9 of
Basrak (2000), the strong mixing condition implies the mixing condition
$\mcr{A}(a_n)$, i.e., there exists a sequence of positive integers $r_n$
such that $r_n\to \infty$, $l_n=[n/r_n]\to \infty$ as $n\to \infty$ and
 \beqnn
\mbf{E}\exp\Big\{-\sum_{k=1}^nf(\mbf{Y}_k/a_n)\Big\}- \Big(\mbf{E}
\exp\Big\{-\sum_{k=1}^{r_n}f(\mbf{Y}_k/a_n)\Big\}\Big)^{l_n} \to0
 \eeqnn
for every $f\in C_0^+(\bar{\mbb{R}}_0)$. In fact, if $\{\mbf{Y}_k\}$ is
strongly mixing with geometric rate, we can choose $r_n=[n^\delta]$ for
any $0<\delta<1$; see Remark~2.3.10 of Basrak (2000). \eremark

Some regularly varying distributions can be found in the SCIR-model.
Recall that $\mbf{P}_x$ denotes the law of the SCIR-model $\{X_t: t\ge
0\}$ defined by (\ref{1.2}) with $X_0=x\ge 0$ and $\mbf{E}_x$ denotes the
corresponding expectation. In the sequel of this section, we assume
$1<\alpha<2$.

\bproposition\label{t3.10} For any $x\ge 0$ we have, as $u\to \infty$,
 \beqnn
\mbf{P}_x(X_t>u)
 \sim
\frac{\sigma^\alpha t}{\alpha\Gamma(-\alpha)} [q_\alpha(t)+p_\alpha(t)x]
u^{-\alpha},
 \eeqnn
where
 \beqlb\label{3.2}
\ar\ar~ p_\alpha(t) = \frac{1}{b(\alpha-1)}[e^{-b t} - e^{-\alpha b t}],
 \quad
q_\alpha(t) = \frac{a}{b}\Big[\frac{1}{\alpha b}(1 - e^{-\alpha b t}) -
p_\alpha(t)\Big].
 \eeqlb
\eproposition

\proof In view of (\ref{1.5}), the extremal behavior of $X_t$ is
determined by a stochastic integral. Then, using Remark~\ref{t3.8}, we
have, as $u\rightarrow \infty$,
 \beqnn
\mbf{P}_x(X_t>u)
 \ar\sim\ar
\mbf{P}_x\Big(\sigma\int_0^te^{-b(t-s)}\sqrt[\alpha]{X_{s-}} dZ_s> u\Big)
\cr
 \ar\sim\ar
\sigma^{\alpha} \mbf{P}_x(Z_t>u) \int_0^te^{-\alpha b(t-s)} \mbf{E}_x(X_s)
ds.
 \eeqnn
Based on (\ref{2.7}), it is easy to compute
 \beqlb\label{3.3}
\mbf{E}_x\Big(\int_0^te^{-\alpha b(t-s)} X_s ds\Big)
 =
q_\alpha(t) + p_\alpha(t)x.
 \eeqlb
By Remark~\ref{t3.7} we have $\mbf{P}_x(Z_t>u)\sim t/\alpha\Gamma(-\alpha)
u^\alpha$. Then the desired result follows. \qed

\bproposition\label{t3.11} For any $K>0$, we have
 \beqnn
\lim_{u\to \infty} \sup_{x\in[0,K]} \Big|u^{\alpha}\mbf{P}_x(X_1>u) -
\frac{\sigma^\alpha}{\alpha\Gamma(-\alpha)}(q_\alpha+p_\alpha x)\Big| = 0
 \eeqnn
and
 \beqnn
\lim_{u\to \infty} \sup_{x\in[0,K]} \Big|u^{\alpha-1} \mbf{E}_x(X_1
1_{\{X_1>u\}}) - \frac{\sigma^\alpha}{(\alpha-1)\Gamma(-\alpha)}
(q_\alpha+p_\alpha x)\Big| = 0.
 \eeqnn
where $p_\alpha = p_\alpha(1)$ and $q_\alpha = q_\alpha(1)$ are defined
by (\ref{3.2}). \eproposition

\proof Let $\{X_t(x)\}$ be the SCIR-model defined by (\ref{1.2}) with
initial value $X_0 = x$. By Theorem~5.5 of Fu and Li (2010), the random
function $x\mapsto X_t(x)$ is increasing, so $x\mapsto \mbf{P}_x(X_t> u)$
is increasing for any $t,u\ge 0$. Then the first convergence holds by
Proposition~\ref{t3.10} and Dini's theorem. The second convergence
follows similarly by Proposition~\ref{t3.2}. \qed

Let us consider a stationary c\`{a}dl\`{a}g realization $\{X_t: t\in
\mbb{R}\}$ of the SCIR-model with one-dimensional marginal distribution
$\mu$ given by (\ref{2.19}). By a modification of the arguments in the
proofs of Theorems~9.31 and~9.32 Li (2011) one can see that, on an
extension of the probability space, there is a compensated Poisson random
measure $\tilde{N}(ds,dz)$ on $\mbb{R}\times (0,\infty)$ with intensity
$ds\nu_\alpha(dz)$ so that (\ref{2.17}) is satisfied for all $t\ge r\in
\mbb{R}$. For any integer $k\in \mbb{Z}$ let
 \beqlb\label{3.4}
\mbf{I}_k = \varepsilon_k\big(1,(1+X_{k-1})^{-1}\big),
 \quad
\mbf{H}_k = X_{k-1}\big(X_{k-1}^{1/\alpha}, \varepsilon_k\big),
 \eeqlb
where
 \beqlb\label{3.5}
\varepsilon_k = \sigma\int_{k-1}^k\int_0^\infty e^{-b(k-s)}
 X_{s-}^{1/\alpha}z \tilde{N}(ds,dz).
 \eeqlb
It is easy to see that the above sequences are stationary. We are going
to prove that the sequences $\{X_k\}$, $\{\mbf{I}_k\}$ and
$\{\mbf{H}_k\}$ are jointly regularly varying.

\bproposition\label{t3.12} Let $\mu$ be the stationary distribution of
the SCIR-model given by (\ref{2.19}). For any $t\ge 0$ we have, as $x\to
\infty$,
 \beqnn
\mu(x,\infty)
 \sim
- \frac{a\sigma^\alpha}{\alpha^2b^2}\Gamma(1-\alpha)^{-1} x^{-\alpha}
 =
\frac{a\sigma^\alpha}{\alpha^3b^2}\Gamma(-\alpha)^{-1} x^{-\alpha}.
 \eeqnn
Consequently, for any $0<r<\alpha$ we have
 \beqnn
\int_0^\infty x^r \mu(dx) = \int_0^\infty \mu(y^{1/r},\infty) dy< \infty.
 \eeqnn
\eproposition

\proof The tail behavior of $X_t$ is closed related with the asymptotics
of its Laplace transform. By (\ref{2.19}), as $\lambda\to 0$,
 \beqnn
L_\mu(\lambda)
 \ar=\ar
1 - \frac{a}{b}\int_0^\lambda \frac{\alpha bdz}{\alpha b+ \sigma^\alpha
z^{\alpha-1}} + O(\lambda^2) \cr
 \ar=\ar
1 - \frac{a}{b}\lambda + \frac{a}{b}\int_0^\lambda \frac{\sigma^\alpha
z^{\alpha-1} dz}{\alpha b+ \sigma^\alpha z^{\alpha-1}} + O(\lambda^2) \cr
 \ar=\ar
1 - \frac{a}{b}\lambda + \frac{a}{b^2}\int_0^\lambda \frac{b\sigma^\alpha
z^{\alpha-1} dz}{\alpha b+ \sigma^\alpha z^{\alpha-1}} + O(\lambda^2) \cr
 \ar=\ar
1 - \frac{a}{b}\lambda + \frac{a}{\alpha b^2}\int_0^\lambda \sigma^\alpha
z^{\alpha-1} dz - \frac{a}{\alpha b^2}\int_0^\lambda
\frac{\sigma^{2\alpha} z^{2(\alpha-1)}dz}{\alpha b+ \sigma^\alpha
z^{\alpha-1}} + O(\lambda^2) \cr
 \ar=\ar
1 - \lambda\int_0^\infty x\mu(dx) + \frac{a\sigma^\alpha}{\alpha^2b^2}
\lambda^\alpha - O(\lambda^{2\alpha-1}) + O(\lambda^2),
 \eeqnn
where we have used (\ref{2.12}) for the last equality. Then the result
follows by Theorem~8.1.6 of Bingham et al.\ (1987). \qed

\blemma\label{t3.13} Let $\tilde{N}(ds,dz)$ be the compensated Poisson
random measure in (\ref{2.17}) and let
 \beqlb\label{3.6}
z(t) = \sigma\int_0^t\int_0^1 e^{-b(t-s)} X_{s-}^{1/\alpha}z
\tilde{N}(ds,dz).
 \eeqlb
Then for any $1\le r< \alpha^2$ and any $T\ge 0$, we have
 \beqnn
\mbf{E}\Big[\sup_{0\le t\le T}|z(t)|^{r}\Big]< \infty.
 \eeqnn
\elemma

\proof We follow an idea in the proof of Lemma~5.5 of Hult and Lindskog
(2007). Let $q = 1/(1-p^{-1})$ for any $p\in (1,\alpha^2/r)$. By the
Burkholder-Davis-Gundy inequality and H\"{o}lder's inequality, we have
 \beqnn
\mbf{E}\Big[\sup_{0\le t\le T}|z(t)|^r\Big]
 \ar\leq\ar
\sigma^r\mbf{E}\Big[\sup_{0\le t\le T}\Big|\int_0^t\int_0^1 e^{bs}
X_{s-}^{1/\alpha}z \tilde{N}(ds,dz)\Big|^r\Big]\cr
 \ar\leq\ar
C_1\mbf{E}\Big[\Big(\int_0^T\int_0^1 e^{2bs}X_{s-}^{2/\alpha}z^2
N(ds,dz)\Big)^{r/2}\Big]\cr
 \ar\leq\ar
C(T)\mbf{E}\Big[ \sup_{0\leq s\leq T}X_s^{r/\alpha} \Big(\int_0^T\int_0^1
z^2 N(ds,dz)\Big)^{r/2}\Big]\cr
 \ar\leq\ar
C(T)\mbf{E}^{1/p} \Big[\sup_{0\leq s\leq T}X_s^{rp/\alpha}\Big]
\mbf{E}^{1/q}\Big[\Big(\int_0^T\int_0^1z^2 N(ds,dz)\Big)^{rq/2}\Big].
 \eeqnn
The first expectation on the right-hand side is finite by
Proposition~\ref{t2.10}. By Theorem~34 of Protter (2005, p.25), the second
expectation is also finite. \qed

\blemma\label{t3.14} Suppose that $\{A_n\}\subset \mcr{F}_0$ is a
sequence of events so that $\mbf{P}(A_n)\to 0$ as $n\to \infty$. Then for
any $x>0$ and $T\ge 0$ we have
 \beqlb\label{3.7}
\lim_{n\to \infty} n\mbf{P}\Big(A_n, \sup_{0\le t\le T} \Big|\int_0^t
e^{-b(t-s)} X_{s-}^{1/\alpha} dZ_s\Big|> a_nx\Big) = 0.
 \eeqlb
\elemma

\proof The right-hand side of (\ref{3.7}) is bound above by
$J_1+J_2+J_3$, where
 \beqnn
J_1 = n\mbf{P}\Big(\sup_{0\le t\le T}\Big|\int_0^t\int_0^1 e^{-b(t-s)}
X_{s-}^{1/\alpha}z \tilde{N}(ds,dz)\Big|> a_nx/3\Big),
 \eeqnn
 \beqnn
J_2 = n\mbf{P}\Big(\sup_{0\le t\le T}\int_0^t e^{-b(t-s)} X_s^{1/\alpha}ds
\int_1^\infty z \nu_\alpha(dz)\Big|> a_nx/3\Big),
 \eeqnn
and
 \beqnn
J_3 \ar:=\ar n\mbf{P}\Big(A_n, \sup_{0\le t\le T}\int_0^t\int_1^\infty
e^{-b(t-s)} X_{s-}^{1/\alpha}z N(ds,dz)\Big|> a_nx/3\Big).
 \eeqnn
Let $z(t)$ be defined by (\ref{3.6}). Then for any $1\le r<\alpha^2$ we
have
 \beqnn
J_1 \le \frac{C_1n}{n^{r/\alpha}x^r}\mbf{E}\Big[\sup_{0\le t\le
T}|z(t)|^r\Big],
 \quad
J_2 \le \frac{C_2n}{n^{r/\alpha}x^r}\mbf{E}\Big[\Big(\int_0^T
X_{s-}^{1/\alpha} ds\Big)^r\Big],
 \eeqnn
where the two expectations are finite by Proposition~\ref{t2.10} and
Lemma~\ref{t3.13}. Then $J_1\to 0$ and $J_2\to 0$ as $n\to \infty$. By
introducing the L\'evy process
 \beqnn
\xi(t) := \int_0^t\int_1^\infty z N(ds,dz), \qquad t\ge 0,
 \eeqnn
for any $K\ge 1$ we have
 \beqnn
J_3 \leq n\mbf{P}(A_n, K^{1/\alpha}\xi(T)> a_nx/6) + n\mbf{P}
\Big(\int_0^T X_{s-}^{1/\alpha} 1_{\{X_{s-}> K\}} d\xi(s)> a_nx/6\Big).
 \eeqnn
By Remark~\ref{t3.7} and the property of independent increments of
$\{\xi(t)\}$ it follows that
 \beqnn
\lim_{n\to \infty}n\mbf{P}(A_n, K^{1/\alpha}\xi(T)> a_nx/6)
 =
\lim_{n\to \infty}n\mbf{P}(A_n)\mbf{P}(K^{1/\alpha}\xi(T)> a_nx/6) = 0.
 \eeqnn
By Remarks~\ref{t3.7} and~\ref{t3.8},
 \beqnn
\lim_{n\to \infty} n\mbf{P}\Big(\int_0^T X_{s-}^{1/\alpha} 1_{\{X_{s-}>
K\}} d\xi(s)> a_nx/6\Big)
 =
C_3x^{-\alpha} \int_0^T \mbf{E}[X_s1_{\{X_s> K\}} ds,
 \eeqnn
which tends to zero as $K\to \infty$. Then we have the desired result.
\qed

\ttheorem\label{t3.15} The sequence $\{X_k\}$ is jointly regular varying
with index $\alpha$. More precisely, as $x\to \infty$ we have
 \beqlb\label{3.8}
\mbf{P}(X_0>x)
 \sim
\frac{a\sigma^\alpha}{\alpha^3b^2}\Gamma(-\alpha)^{-1} x^{-\alpha}
 \eeqlb
and, for any integer $k\ge 1$,
 \beqlb\label{3.9}
\mbf{P}\big(X_0^{-1}(X_0,\cdots,X_k)\in \cdot\big|X_0> x\big)
 \overset{\rm w}{\longrightarrow}
\delta_{(1,e^{-b},\cdots,e^{-bk})}(\cdot).
 \eeqlb
\etheorem

\proof By Proposition~\ref{t3.12} we have the asymptotics (\ref{3.8}). It
suffices to show (\ref{3.9}) holds when $x\to \infty$ along the sequence
$a_n := n^{1/\alpha}$. Let $\mbf{X}=(X_0,X_1,\cdots,X_k)$ and
$\bar{\mbf{X}} = (X_0, X_0e^{-b}, \cdots, X_0e^{-bk})$. Let $z(t)$ be
defined by (\ref{3.6}). For any $\delta>0$, we can use (\ref{1.5}) to see
 \beqnn
\ar\ar\mbf{P}\big(\|\mbf{X}-\bar{\mbf{X}}\|> a_n\delta\big|X_0> a_n\big)
\cr
 \ar\ar\qquad
\le \mbf{P}\Big(ak + \max_{1\le j\le k}\Big|\int_0^j e^{-b(j-s)}
X_{s-}^{1/\alpha} dZ_s\Big|> a_n\delta\Big|X_0> a_n\Big) \cr
 \ar\ar\qquad
\le Cn\mbf{P}\Big(X_0> a_n, \max_{1\le j\le k}\Big|\int_0^j e^{-b(j-s)}
X_{s-}^{1/\alpha} dZ_s\Big|> a_n\delta - ak\Big) \cr
 \ar\ar\qquad
\le Cn\sum_{j=1}^k \mbf{P}\Big(X_0> a_n, \Big|\int_0^j e^{-b(j-s)}
X_{s-}^{1/\alpha} dZ_s\Big|> (a_n\delta-ak)/k\Big).
 \eeqnn
By Lemma~\ref{t3.14} it is easy to see the right-hand tends to zero as
$n\to \infty$. By Theorem~\ref{t3.3} we conclude that, as $n\to \infty$,
 \beqnn
\mbf{P}(\mbf{X}/a_n\in\cdot|X_0> a_n)
 \ar\sim\ar
\mbf{P}(\bar{\mbf{X}}/a_n\in\cdot|X_0> a_n) \cr
 \ar\overset{\rm v}{\longrightarrow}\ar
\alpha\int_1^\infty 1_{\{z(1,e^{-b},\cdots, e^{-bk})
\in\cdot\}}\frac{dz}{z^{\alpha+1}}.
 \eeqnn
Then (\ref{3.9}) follows by the continuous mapping theorem. By
Theorem~\ref{t3.5}, the sequence $\{X_k\}$ is jointly regular varying
with index $\alpha$. \qed

\ttheorem\label{t3.16} The sequence $\{\mbf{I}_k\}$ defined by
(\ref{3.4}) is jointly regular varying with index $\alpha$. More
precisely, as $x\to \infty$ we have
 \beqlb\label{3.10}
\mbf{P}(\|\mbf{I}_1\|>x)
 \sim
\frac{\mbf{E}(G)}{\alpha^2\Gamma(-\alpha)} x^{-\alpha}
 =
\frac{a\sigma^\alpha(1-e^{-\alpha b})}{\alpha^3b^2\Gamma(-\alpha)}
 x^{-\alpha}
 \eeqlb
and, for any integer $k\ge 1$,
 \beqlb\label{3.11}
\ar\ar\mbf{P}\big(\|\mbf{I}_1\|^{-1}(\mbf{I}_1,\cdots, \mbf{I}_k)\in
\cdot \big|\|\mbf{I}_1\|> x\big) \ccr
 \ar\ar\qquad
\overset{\rm w}{\longrightarrow} \mbf{E}(G)^{-1}\mbf{E}\big[G;
((1,(1+X_0)^{-1}),\mbf{0},\cdots,\mbf{0})\in \cdot\big],
 \eeqlb
where
 \beqlb\label{3.12}
G=\sigma^\alpha\int_0^1e^{-\alpha b(1-t)}X_tdt.
 \eeqlb
\etheorem

\proof Clearly, it suffices to show that (\ref{3.10}) and (\ref{3.11})
hold when $x\to \infty$ along the sequence $a_n := n^{1/\alpha}$. For
$0\le s\le k$, let $\mbf{Z}_s=(Z_s,Z_s)$ and $\mbf{\Phi}_s =
(\phi_1(s),\phi_2(s))$, where
 \beqnn
\phi_1(s)=\sum_{j=1}^k1_{(j-1,j]}(s)e^{-b(j-s)} \sigma X_{s-}^{1/\alpha},
 ~
\phi_2(s)=\sum_{j=1}^k1_{(j-1,j]}(s)\frac{e^{-b(j-s)} \sigma
X_{s-}^{1/\alpha}}{1+X_{j-1}}.
 \eeqnn
We consider the process
 \beqnn
(\mbf{\Phi}\cdot\mbf{Z})_t
 :=
\Big(\int_0^t\phi_1(s)dZ_s,\int_0^t\phi_2(s)dZ_s\Big).
 \eeqnn
Choose some $\delta\in (\alpha,\alpha^2)$. By Proposition~\ref{t3.12}, we
have $\mbf{E}[X_0^{\delta/\alpha}]<\infty$. Then Proposition~\ref{t2.10}
implies that
 \beqnn
\mbf{E}\Big[\sup_{0\le s\le k}\|\mbf{\Phi}_s\|^{\delta}\Big]
 \le
\mbf{E}\Big[\sup_{0\le s\le k}\sigma^\delta X_s^{\delta/\alpha}\Big]
 \le
\sigma^{\delta}C(k)[1+\mbf{E}(X_0^{\delta/\alpha})]< \infty.
 \eeqnn
By Theorem~3.4 in Hult and Lindskog (2007) and Remark~\ref{t3.7}, as
$n\to \infty$,
 \beqlb\label{3.13}
\ar\ar~ n\mbf{P}(a_n^{-1}(\mbf{\Phi}\cdot \mbf{Z})\in\cdot)
 \overset{\rm \hat{w}}{\longrightarrow}
Q(\cdot) := k\mbf{E}[\nu_\alpha\{x\in\mbb{R}_+: x\mbf{\Phi}_\tau
1_{[\tau,k]}\in\cdot\}],
 \eeqlb
where $\nu_\alpha$ is defined by (\ref{1.3}) and $\tau$ is uniformly
distributed on $[0,k]$ and independent of $\bf\Phi$. In view of
(\ref{3.13}), we have
 \beqnn
n\mbf{P}(a_n^{-1}(\mbf{\Phi}\cdot\mbf{Z})_1\in\cdot)
 \overset{\rm v}{\longrightarrow}
Q(\{\mbf{y}\in \mbb{D}^2[0,k]: \mbf{y}_1\in\cdot\}).
 \eeqnn
By Definition~\ref{t3.1} it follows that
 \beqnn
Q(\{\mbf{y}\in \mbb{D}^2[0,k]: \|\mbf{y}_1\|>r\})
 =
r^{-\alpha}Q(\{\mbf{y}\in \mbb{D}^2[0,k]: \|\mbf{y}_1\|>1\}).
 \eeqnn
Now define the functions $h_0,h_1,h_2:\mbb{D}^2[0,k]\to \mbb{R}^{2k}$ by
 \beqnn
\ar\ar h_0(\mbf{y}) = (\mbf{y}_1, \mbf{y}_2-\mbf{y}_1,\cdots,
\mbf{y}_k-\mbf{y}_{k-1}), \ccr
 \ar\ar
h_1(\mbf{y}) = 1_{\{\|\mbf{y}_1\|>1\}},
 ~~
h_2(\mbf{y}) = h_0(\mbf{y})h_1(\mbf{y}).
 \eeqnn
Let ${\rm Disc}(h_i)$ be the set of discontinuities of $h_i$ ($i=0,1,2$).
By (\ref{3.13}) it is easy to see that $Q({\rm Disc}(h_0)) = Q({\rm
Disc}(h_1)) = 0$, so $Q({\rm Disc}(h_2))=0$. Moreover, for any $B\in
\mcr{B}(\mbb{R}^{2k})$ bounded away from $\bf 0$ the set $h^{-1}_2(B)\in
\mcr{B}(\mbb{D}^2[0,k])$ is bounded away from $\bf 0$. Applying the
continuous mapping theorem, we obtain as $n\to \infty$,
 \beqnn
n\mbf{P}\big(\|\mbf{I}_1\|> a_n,\
a_n^{-1}(\mbf{I}_1,\mbf{I}_2,\cdots,\mbf{I}_k)\in \cdot\big)
 \overset{\rm v}{\longrightarrow}
Q\circ h_2^{-1}(\cdot)
 \eeqnn
on $\mbb{R}^{2k}\setminus \{\mbf{0}\}$, where
 \beqnn
Q\circ h_2^{-1}(\cdot)
 \ar=\ar
k\mbf{E}\Big[\nu_\alpha\big\{x\in\mbb{R}_+:\|x\mbf{\Phi}_\tau
1_{[\tau,k]}(1)\|>1, \cr
 \ar\ar\qquad\qquad\qquad\qquad\qquad\qquad
x\big(\mbf{\Phi}_\tau 1_{[\tau,k]}(1),\mbf{0},\cdots,\mbf{0}\big)\in
\cdot\big\}\Big] \cr
 \ar=\ar
k\mbf{E}\Big[\nu_\alpha\big\{x\in\mbb{R}_+:\|x\mbf{\Phi}_\tau\|>1,
\tau\le 1, \cr
 \ar\ar\qquad\qquad\qquad\qquad\qquad\qquad
x\big(\mbf{\Phi}_\tau 1_{[\tau,k]}(1),\mbf{0},\cdots,\mbf{0}\big)\in
\cdot\big\}\Big] \cr
 \ar=\ar
k\mbf{E}\Big[\nu_\alpha\big\{x\in\mbb{R}_+:x\|\mbf{\Phi}_\tau\|>1, x
\big(\mbf{\Phi}_\tau,\mbf{0},\cdots,\mbf{0}\big)\in
\cdot\big\}1_{[0,1]}(\tau)\Big].
 \eeqnn
Let $E=\{(\mbf{y}_1,\cdots,\mbf{y}_k)\in\mbb{R}^{2k}:\|\mbf{y}_1\|>0\}$.
Define the injection $f: E\to(0,\infty)\times E$ by
 \beqnn
f(\mbf{y}_1,\cdots,\mbf{y}_k)
 =
\big(\|\mbf{y}_1\|,\mbf{y}_1/\|\mbf{y}_1\|,\cdots,\mbf{y}_k/\|
\mbf{y}_1\|\big).
 \eeqnn
Then we have as $n\to \infty$,
 \beqnn
n\mbf{P}\big(\|\mbf{I}_1\|> a_n,\
\|\mbf{I}_1\|^{-1}(\mbf{I}_1,\cdots,\mbf{I}_k) \in \cdot\big)
 \overset{\rm w}{\longrightarrow}
Q\circ h_2^{-1}\circ f^{-1}(\cdot)
 \eeqnn
on $E$, where
 \beqnn
\ar\ar Q\circ h_2^{-1}\circ f^{-1}(\cdot) \ccr
 \ar\ar\qquad
= k\mbf{E}\Big[\nu_\alpha\big\{x\in \mbb{R}_+: x\|\mbf{\Phi}_\tau\|> 1,
\big(\mbf{\Phi}_\tau/\|\mbf{\Phi}_\tau\|,\mbf{0},\cdots,\mbf{0}\big)\in
\cdot\big\}1_{[0,1]}(\tau)\Big]\cr
 \ar\ar\qquad
= k\mbf{E}\Big[\nu_\alpha\big\{x\in \mbb{R}_+: x\|\mbf{\Phi}_\tau\|>
1\big\} 1_{[0,1]}(\tau), \big(\mbf{\Phi}_\tau/\|\mbf{\Phi}_\tau\|,
\mbf{0}, \cdots, \mbf{0}\big)\in \cdot\Big]\cr
 \ar\ar\qquad
= \frac{k}{\alpha\Gamma(-\alpha)}
\mbf{E}\Big[\int_{\|\mbf{\Phi}_\tau\|^{-1}}^\infty
\frac{dz}{z^{\alpha+1}} 1_{[0,1]}(\tau),\big((1,(1+X_0)^{-1}),
\mbf{0},\cdots,\mbf{0}\big)\in \cdot\Big]\cr
 \ar\ar\qquad
= \frac{k}{\alpha^2\Gamma(-\alpha)}
\mbf{E}\big[\|\mbf{\Phi}_\tau\|^{\alpha} 1_{[0,1]}(\tau),
\big((1,(1+X_0)^{-1}), \mbf{0},\cdots,\mbf{0}\big)\in \cdot\big]\cr
 \ar\ar\qquad
= \frac{1}{\alpha^2\Gamma(-\alpha)}\int_0^1
\mbf{E}\big[\|\mbf{\Phi}_s\|^{\alpha},
\big((1,(1+X_0)^{-1}),\mbf{0},\cdots, \mbf{0}\big)\in \cdot\big]ds \cr
 \ar\ar\qquad
= \frac{1}{\alpha^2\Gamma(-\alpha)} \mbf{E}\big[G, \big((1,(1+X_0)^{-1}),
\mbf{0},\cdots,\mbf{0}\big)\in \cdot\big].
 \eeqnn
In particular, as $n\to \infty$ we have
 \beqnn
n\mbf{P}(\|\mbf{I}_1\|> a_n)\longrightarrow \frac{\mbf{E}(G)}
{\alpha^2\Gamma(-\alpha)}.
 \eeqnn
By (\ref{2.12}), we have $\mbf{E}(X_t)=a/b$. It follows that
 \beqnn
\mbf{E}(G)
 =
\frac{a\sigma^\alpha}{b}\int_0^1e^{-\alpha b(1-t)} dt
 =
\frac{a\sigma^\alpha}{\alpha b^2}(1-e^{-\alpha b}).
 \eeqnn
Then we have (\ref{3.10}) and (\ref{3.11}). By Theorem~\ref{t3.5}, the
sequence $\{\mbf{I}_k\}$ is jointly regular varying with index $\alpha$.
\qed

\bremark\label{t3.17} For $k=1,2,\cdots$ define
 \beqlb\label{3.14}
V_k = \sigma \int_{k-1}^ke^{-b(k-s)}e^{-b(s-k+1)/\alpha} dZ_s.
 \eeqlb
Then the sequence $\{V_k\}$ is i.i.d.\ with the same distribution as
 \beqnn
\sigma\Big(\frac{e^{-b} - e^{-\alpha b}} {(\alpha-1)b}\Big)^{1/\alpha}
Z_1,
 \eeqnn
which is regularly varying with index $\alpha$. \eremark

\blemma\label{t3.18} Let $V_k$ be defined by (\ref{3.14}) and let
$\bar{\mbf{H}}_k = X_{k-1}^{(\alpha + 1)/\alpha} (1,V_k)$. Then for any
$0< r< \alpha^3/(\alpha^2+1)$, we have
 \beqnn
\mbf{E}\big[\|\mbf{H}_k-\bar{\mbf{H}}_k\|^r\big]< \infty.
 \eeqnn
\elemma

\proof Since $0<r< \alpha^3(\alpha^2+1)^{-1}< \alpha$, by
Lemma~\ref{t2.8} and H\"older's inequality,
 \beqnn
\mbf{E}\big[\|\mbf{H}_k-\bar{\mbf{H}}_k\|^r\big]
 \ar=\ar
\mbf{E}\Big[\Big|\int_{k-1}^k\sigma
e^{-b(k-s)}X_{k-1}\Big(X_{s-}^{1/\alpha} - X_{k-1}^{1/\alpha}
e^{-b(s-k+1)/\alpha}\Big)dZ_s\Big|^r\Big] \cr
 \ar\le\ar
C\mbf{E}\Big[\Big(\int_{k-1}^k e^{-\alpha b(k-s)}X_{k-1}^{\alpha}\big|X_s
-X_{k-1}e^{-b(s-k+1)}\big|ds\Big)^{r/\alpha}\Big] \cr
 \ar\le\ar
C\mbf{E}\Big\{X_{k-1}^r \Big[\mbf{E}_{X_{k-1}}\Big(\int_0^1 |X_s - X_0
e^{-bs}|ds\Big)\Big]^{r/\alpha}\Big\}\cr
 \ar\le\ar
C\mbf{E}\Big\{X_{k-1}^r \Big[\int_0^1 \mbf{E}_{X_{k-1}}\big(|X_s - X_0
e^{-bs}|\big)ds\Big]^{r/\alpha}\Big\}\cr
 \ar\le\ar
C\mbf{E}\Big\{X_{k-1}^r\Big[\int_0^1 \mbf{E}_{X_{k-1}} \Big(\Big|\int_0^s
e^{-b(s-u)} X_{u-}^{1/\alpha} dZ_u\Big|\Big) ds\Big]^{r/\alpha}\Big\} \cr
 \ar\le\ar
C\mbf{E}\Big\{X_{k-1}^r\Big[\int_0^1 \mbf{E}_{X_{k-1}} \Big(\Big|\int_0^s
e^{-b(s-u)} X_{u-}^{1/\alpha} dZ_u\Big|^r\Big) ds\Big]^{1/\alpha}\Big\}
\cr
 \ar\le\ar
C\mbf{E}\big[X_{k-1}^r(1 + X_{k-1}^{r/\alpha^2})\big],
 \eeqnn
which is finite by Proposition~\ref{t2.10}. \qed

\ttheorem\label{t3.19} The sequence $\{\mbf{H}_k\}$ defined by
(\ref{3.4}) is jointly regular varying with index $\alpha^2/(\alpha+1)$.
Let $V_1$ be defined by (\ref{3.14}). Then we have, as $x\to \infty$,
 \beqlb\label{3.15}
\mbf{P}(\|\mbf{H}_1\|>x)
 \sim
\mbf{E}\big[1\vee |V_1|^{\alpha^2/(\alpha+1)}\big]
\frac{a\sigma^\alpha}{\alpha^3b^2\Gamma(-\alpha)} x^{-\alpha^2/(\alpha+1)}
 \eeqlb
and, for any integer $k\ge 1$,
 \beqlb\label{3.16}
\ar\ar\mbf{P}\big(\|\mbf{H}_1\|^{-1}(\mbf{H}_1,\cdots, \mbf{H}_k)\in
\cdot \big|\|\mbf{H}_1\|>x\big) \cr
 \ar\ar\qquad
\overset{\rm w}{\longrightarrow} \frac{\mbf{E}\big[1\vee
|V_1|^{\alpha^2/(\alpha+1)}; (\mbf{\Theta}_1, \cdots,
\mbf{\Theta}_k)\in\cdot\big]} {\mbf{E}\big[1\vee
|V_1|^{\alpha^2/(\alpha+1)}\big]},
 \eeqlb
where
 \beqnn
\mbf{\Theta}_j = e^{-b(j-1)(\alpha+1)/\alpha}(1\vee |V_1|)^{-1}(1,V_j).
 \eeqnn
 \etheorem

\proof We only need to show (\ref{3.15}) and (\ref{3.16}) hold when $x\to
\infty$ along the sequence $c_n := n^{(\alpha+1)/\alpha^2}$. By
Proposition~\ref{t3.12} it follows that $X_0^{(\alpha+1)/\alpha}$ is
regularly varying with the index $\alpha^2/(\alpha+1)$. More precisely,
as $n\to \infty$,
 \beqlb\label{3.17}
\mbf{P}\big(X_0^{(\alpha+1)/\alpha}> x\big)
 \sim
\frac{a\sigma^\alpha}{\alpha^3b^2\Gamma(-\alpha)}
x^{-\alpha^2/(\alpha+1)}.
 \eeqlb
By Remark~\ref{t3.17}, we have $\mbf{E}[|V_k|^r]< \infty$ for any
$0<r<\alpha$. Note that $\|\bar{\mbf{H}}_k\| =
X_{k-1}^{(\alpha+1)/\alpha} (1\vee V_k)$. By (\ref{3.17}) and Breiman's
Lemma, as $n\to \infty$,
 \beqlb\label{3.18}
\mbf{P}(\|\bar{\mbf{H}}_1\|>x)
 \sim
\mbf{E}\big[1\vee |V_1|^{\alpha^2/(\alpha+1)}\big]
\frac{a\sigma^\alpha}{\alpha^3b^2\Gamma(-\alpha)}
x^{-\alpha^2/(\alpha+1)};
 \eeqlb
see, e.g., Resnick (1987, p.231). Then $\|\bar{\mbf{H}}_1\|$ is regularly
varying with the index $\alpha^2/(\alpha+1)$. By Lemma~\ref{t3.18}, for
any $0<r< \alpha^3/(\alpha^2+1)$, we have $\mbf{E}[\|\mbf{H}_1 -
\bar{\mbf{H}}_1\|^r]< \infty$. Now let $\mbf{H}=(\mbf{H}_1, \cdots,
\mbf{H}_k)$ and $\bar{\mbf{H}} = (\bar{\mbf{H}}_1,\cdots,
\bar{\mbf{H}}_k)$. By Markov's inequality and (\ref{3.18}) we have, as
$x\to \infty$,
 \beqnn
\frac{\mbf{P}(\|\mbf{H}_1-\bar{\mbf{H}}_1\|>x)}
{\mbf{P}(\|\bar{\mbf{H}}_1\|>x)}
 \ar\le\ar
\frac{\mbf{P}(\|\mbf{H}-\bar{\mbf{H}}\|>x)}
{\mbf{P}(\|\bar{\mbf{H}}_1\|>x)} \cr
 \ar\le\ar
\frac{x^{-r}}{\mbf{P}(\|\bar{\mbf{H}}_1\|>x)}\sum_{j=1}^k
\mbf{E}\big[\|\mbf{H}_j-\bar{\mbf{H}}_j\|^r\big]\to 0.
 \eeqnn
As in the proof of Lemma~3.12 in Jessen and Mikosch (2006), we obtain
$(\ref{3.15})$ from (\ref{3.18}). From the above relation we also have,
as $n\to \infty$,
 \beqlb\label{3.19}
\mbf{P}(c_n^{-1}\mbf{H}\in\cdot|\|\mbf{H}_1\|>c_n)
 \sim
\mbf{P}(c_n^{-1}\bar{\mbf{H}}\in\cdot|\|\bar{\mbf{H}}_1\|> c_n).
 \eeqlb
Let $\tilde{\mbf{H}}_k = (X_0 e^{-b(k-1)})^{(\alpha+1)/\alpha} (1,V_k)$.
For $\delta> 0$ and $K> 1$ we have
 \beqnn
\ar\ar\mbf{P}\big(\|\tilde{\mbf{H}}_k-\bar{\mbf{H}}_k\|> c_n\delta \big|
X_0> a_n\big) \cr
 \ar\ar\qquad
\le \mbf{P}\big(\|\tilde{\mbf{H}}_k-\bar{\mbf{H}}_k\|> c_n\delta, |V_k|\le
K \big| X_0> a_n\big) \cr
 \ar\ar\qquad\quad
+\, \mbf{P}\big(\|\tilde{\mbf{H}}_k-\bar{\mbf{H}}_k\|> c_n\delta, |V_k|> K
\big| X_0> a_n\big).
 \eeqnn
Let $J_1$ and $J_2$ denote the two terms on the right-hand side. Then
 \beqnn
J_1\le \mbf{P}\big(K|(X_0e^{-b(k-1)})^{(\alpha+1)/\alpha} -
X_{k-1}^{(\alpha+1)/\alpha}|> c_n\delta\big| X_0> a_n\big).
 \eeqnn
By Theorem~\ref{t3.15} and the continuous mapping theorem, we have
$J_1\to 0$ as $n\to \infty$. Since $X_0$ is independent of $V_k$, we have
 \beqnn
J_2\le \mbf{P}(|V_k|>K|X_0> a_n) = \mbf{P}(|V_k|>K),
 \eeqnn
which tends to zero as $K\to \infty$. Then the regular variation property
of $X_0$ implies, for any $\zeta>0$,
 \beqlb\label{3.20}
\lim_{x\to \infty}\mbf{P}(\|\tilde{\mbf{H}}_k-\bar{\mbf{H}}_k\|>
x^{(\alpha+1)/\alpha}\zeta\big|X_0>x)=0.
 \eeqlb
See, e.g., Resnick (1987, p.14) for a similar method. By (\ref{3.8}) we
have $\mbf{P}(K^{\alpha/(\alpha+1)}X_0> a_n)\sim h_1n^{-1}$ as $n\to
\infty$ for some constant $h_1>0$. By (\ref{3.20}) and the multiplicative
formula it is easy to see
 \beqnn
\lim_{n\to \infty}n\mbf{P}(\|\tilde{\mbf{H}}_k-\bar{\mbf{H}}_k\|>
c_n\delta, K^{\alpha/(\alpha+1)}X_0> a_n) = 0.
 \eeqnn
It follows that
 \beqnn
\ar\ar \limsup_{n\to \infty}n\mbf{P}(\|\tilde{\mbf{H}}_k -
\bar{\mbf{H}}_k\|> c_n\delta, \|\bar{\mbf{H}}_1\|> c_n)\cr
 \ar\ar\qquad
= \limsup_{n\to \infty}n\mbf{P}(\|\tilde{\mbf{H}}_k - \bar{\mbf{H}}_k\|>
c_n\delta, X_0^{(\alpha+1)/\alpha} (1\vee V_1)> c_n) \cr
 \ar\ar\qquad
\le \lim_{n\to \infty}n\mbf{P}(X_0^{(\alpha+1)/\alpha}|V_1|
1_{\{|V_1|>K\}}>c_n) \cr
 \ar\ar\qquad
= \lim_{n\to \infty} n\mbf{E}[|V_1|^{\alpha^2/(\alpha+1)}1_{\{|V_1|>K\}}]
\mbf{P}(X_0^{(\alpha+1)/\alpha}> c_n) \cr
 \ar\ar\qquad
= C_1\mbf{E}[|V_1|^{\alpha^2/(\alpha+1)}1_{\{|V_1|>K\}}],
 \eeqnn
where we have used Breiman's Lemma again for the second equality. The
right hand side goes to zero as $K\to \infty$. But, by (\ref{3.18}) there
is a constant $h_2> 0$ so that $\mbf{P}(\|\bar{\mbf{H}}_1\|> c_n)\sim
h_2n^{-1}$ as $n\to \infty$. Then
 \beqnn
\lim_{n\to \infty} \mbf{P}(\|\tilde{\mbf{H}}_k- \bar{\mbf{H}}_k\|>
c_n\delta| \|\bar{\mbf{H}}_1\|>c_n)=0.
 \eeqnn
Let $\tilde{\mbf{H}}=(\tilde{\mbf{H}}_1,\cdots,\tilde{\mbf{H}}_k)$. We
have
 \beqnn
\mbf{P}(\|\tilde{\mbf{H}}-\bar{\mbf{H}}\|>c_n\delta|
\|\tilde{\mbf{H}}_1\|>c_n) \le\sum_{j=1}^k\mbf{P}(\|\tilde{\mbf{H}}_j
-\bar{\mbf{H}}_j\|> c_n\delta|\|\tilde{\mbf{H}}_1\|>c_n) \to0.
 \eeqnn
Since $\tilde{\mbf{H}}_1=\bar{\mbf{H}}_1$, by the above relation, we have
as $n\to \infty$,
 \beqlb\label{3.21}
\mbf{P}(c_n^{-1}\bar{\mbf{H}}\in\cdot|\|\bar{\mbf{H}}_1\|> c_n)
 \ar\sim\ar
\mbf{P}(c_n^{-1}\tilde{\mbf{H}}\in\cdot|\|\tilde{\mbf{H}}_1\|> c_n) \cr
 \ar\sim\ar
h_2^{-1}n\mbf{P}(|\|\tilde{\mbf{H}}_1\|>c_n, c_n^{-1}
\tilde{\mbf{H}}\in\cdot).
 \eeqlb
By Proposition~\ref{t3.12} one can see, as $n\to \infty$,
 \beqnn
n \mbf{P}(c_n^{-1}X_0^{(\alpha+1)/\alpha}\in\cdot)
 \overset{\rm v}{\longrightarrow}
\nu(\cdot) := C_2\int_0^\infty \frac{1_{\{u\in\cdot\}}du} {u^{1 +
\alpha^2/(\alpha+1)}}.
 \eeqnn
Note that $X_0$ is independent of ${V_k}$ for $k\ge 1$. By the extended
Breiman's Lemma, as $n\to \infty$,
 \beqnn
n\mbf{P}(c_n^{-1}{\bf\tilde{H}}\in\cdot)
 \overset{\rm v}{\longrightarrow}
\mbf{E}\big[\nu\{u: u(\mbf{\Xi}_1,\cdots, \mbf{\Xi}_k)\in \cdot\}\big].
 \eeqnn
where $\mbf{\Xi}_j = e^{-b(j-1)(\alpha+1)/\alpha}(1,V_j)$; see
Theorem~3.1 of Hult and Lindskog (2007). By (\ref{3.19}) and
(\ref{3.21}), we have
 \beqnn
\ar\ar\mbf{P}(c_n^{-1}\mbf{H}\in\cdot|\|\mbf{H}_1\|> c_n) \cr
 \ar\ar\qquad
\overset{\rm v}{\longrightarrow} C_3\int_0^\infty
\mbf{E}\big[\nu\{u>(1\vee V_1)^{-1}: u(\mbf{\Xi}_1, \cdots,
\mbf{\Xi}_k)\in\cdot\}\big] \frac{du}{u^{1 + \alpha^2/(\alpha+1)}}.
 \eeqnn
Then (\ref{3.16}) follows by an application of the continuous mapping
theorem. By Theorem~\ref{t3.5}, the sequence $\{\mbf{H}_k\}$ is jointly
regular varying with index $\alpha$. \qed


\section{Point processes and partial sums}

\setcounter{equation}{0}

In this section, we will first prove some limit theorems on the point
processes associated with the stationary sequences $\{\mbf{I}_k\}$ and
$\{\mbf{H}_k\}$ defined by (\ref{3.4}). From the limit theorems, we
derive the limits of suitably normalized partial sums for those
sequences. The techniques have been developed extensively by Basrak and
Segers (2009), Davis and Hsing (1995), Davis and Mikosch (1998) among
others. For i.i.d.\ random variables, the idea goes back to Davis (1983)
and LePage et al.\ (1981). Throughout this section, we assume
$1<\alpha<2$. Let $a_n=n^{1/\alpha}$ and $c_n = n^{(\alpha+1)/\alpha^2} =
a_n^{(\alpha+1)/\alpha}$ for $n\ge 1$.

\blemma\label{t4.1} Let $r_n=[n^\delta]$ for any $0<\delta<1$. Then for
any $x>0$ we have
 \beqlb\label{4.1}
\lim_{m\to \infty}\limsup_{n\to \infty}n\sum_{k=m}^{r_n}
\mbf{P}(|\varepsilon_k|> a_nx,|\varepsilon_1|> a_nx)=0.
 \eeqlb
\elemma

\proof Let $\mbf{E}_x$ denote the expectation of $\{X_t: t\ge 0\}$ given
$X_0=x$. Take a constant $r\in (\delta,1)$. For $k\ge 2$, we can use
Markov's inequality and Proposition~\ref{t2.9} to see
 \beqnn
\ar\ar\mbf{P}(|\varepsilon_k|> a_nx,|\varepsilon_1|> a_nx) \cr
 \ar\ar\qquad
= \frac{\sigma^{r\alpha}}{(a_nx)^{r\alpha}}
\mbf{E}\Big[1_{\{|\varepsilon_1|> a_nx\}} \Big|\int_{k-1}^k e^{-b(k-s)}
X_{s-}^{1/\alpha} dZ_s\Big|^{r\alpha}\Big]\cr
 \ar\ar\qquad
= \frac{\sigma^{r\alpha}} {(a_nx)^{r\alpha}}
\mbf{E}\Big[1_{\{|\varepsilon_1|> a_nx\}} \mbf{E}_{X_{k-1}}
\Big(\Big|\int_0^1 e^{-b(1-s)} X_{s-}^{1/\alpha}
dZ_s\Big|^{r\alpha}\Big)\Big]\cr
 \ar\ar\qquad
\le \frac{C_1\sigma^{r\alpha}} {(a_nx)^{r\alpha}}
\mbf{E}\big[1_{\{|\varepsilon_1|> a_nx\}} (1+X_{k-1}^r)\big]\cr
 \ar\ar\qquad
\le \frac{C_2\sigma^{r\alpha}} {(a_nx)^{r\alpha}}
\mbf{E}\big[1_{\{|\varepsilon_1|> a_nx\}} (1+X_1^r e^{-rb(k-2)/\alpha}
)\big],
 \eeqnn
where the last inequality follows from Proposition~\ref{t2.10}. In view
of (\ref{1.6}), we have
 \beqnn
n\sum_{k=m}^{r_n} \mbf{P}(|\varepsilon_k|> a_nx,|\varepsilon_1|> a_nx)\le
J_1+J_2+J_3,
 \eeqnn
where
 \beqnn
J_1\ar=\ar\frac{C_3\sigma^{r\alpha}n}{(a_nx)^{r\alpha}} (r_n-m+1)
\mbf{P}\big(|\varepsilon_1|> a_nx\big), \cr
 J_2\ar=\ar\frac{C_4\sigma^{r\alpha}n}{(a_nx)^{r\alpha}}
\sum_{k=m}^\infty e^{-rb(k-2)/\alpha} \mbf{E}\big(|\varepsilon_1|^r
1_{\{|\varepsilon_1|> a_nx\}}\big), \cr
 J_3\ar=\ar\frac{C_5\sigma^{r\alpha}n}{(a_nx)^{r\alpha}}
\sum_{k=m}^\infty e^{-rb(k-2)/\alpha} \mbf{E}\big(X_0^r
1_{\{|\varepsilon_1|> a_nx\}}\big).
 \eeqnn
By Theorem~\ref{t3.16}, we have $\mbf{P}(|\varepsilon_1|> a_nx) =
O(n^{-1})$. It follows that $J_1 = O(n^{\delta-r})$ as $n\to \infty$. By
Proposition~\ref{t3.2} one can see, as $n\to \infty$,
 \beqnn
\mbf{E}\big(|\varepsilon_1|^r1_{\{|\varepsilon_1|> a_nx\}}\big)
 \ar\sim\ar
\frac{\alpha}{\alpha-r}(a_nx)^r \mbf{P}\big(|\varepsilon_1|> a_nx\big)
\cr
 \ar\sim\ar
\frac{\alpha x^r}{\alpha-r}n^{r/\alpha} \mbf{P}\big(|\varepsilon_1|>
a_nx\big).
 \eeqnn
Thus we have $J_2 = O(n^{r/\alpha-r})$ as $n\to \infty$. By Markov's
inequality and Proposition~\ref{t2.9},
 \beqnn
\mbf{E}\big(X_0^r1_{\{|\varepsilon_1|> a_nx\}}\big)
 \ar\le\ar
\frac{1}{(a_nx)^{\alpha(1-r)}} \mbf{E}\big(X_0^r
|\varepsilon_1|^{\alpha(1-r)}\big) \cr
 \ar=\ar
\frac{1}{(a_nx)^{\alpha(1-r)}}\mbf{E}\big\{X_0^r[\mbf{E}_{X_0}
(|\varepsilon_1|^{\alpha(1-r)})]\big\} \cr
 \ar\le\ar
\frac{C_6}{(a_nx)^{\alpha(1-r)}}\mbf{E}(X_0^r+X_0).
 \eeqnn
It follows that
 \beqnn
\lim_{m\to \infty} \limsup_{n\to \infty}J_3\le \lim_{m\to \infty}
C_7\sum_{k=m}^{\infty} e^{-r\alpha b(k-1)} =0.
 \eeqnn
Then we have (\ref{4.1}). \qed

\ttheorem\label{t4.2} Let $G$ be defined by (\ref{3.12}). Then we have,
as $n\to \infty$,
 \beqlb\label{4.2}
\eta_n := \sum_{k=1}^n\delta_{a_n^{-1}\mbf{I}_k}\overset{\rm
d}{\longrightarrow} \eta \quad\mbox{on}\quad M(\bar{\mbb{R}}_0^2),
 \eeqlb
where $\eta$ is a point process on $\bar{\mbb{R}}_0^2$ with the Laplace
functional $\mbf{E}[e^{-\eta(f)}]$, $f\in C_0^+(\bar{\mbb{R}}_0^2)$ given
by
 \beqlb\label{4.3}
\ar\ar~ \exp\bigg\{-\frac{1}{\alpha\Gamma(-\alpha)}\int_0^\infty
\mbf{E}\Big[\Big(1 - \exp\Big\{-f\Big(y,\frac{y}{1+X_0}\Big)\Big\}\Big)
G\Big]\frac{dy}{y^{\alpha+1}}\bigg\}.
 \eeqlb
\etheorem

\proof By Theorem~\ref{t2.7}, the process $\{X_t\}$ is strongly mixing
with geometric rate. From (\ref{1.6}) and (\ref{3.4}) we see $\mbf{I}_k$
is measurable with respect to $\sigma(X_{k-1},X_k)$. Then $\{\mbf{I}_k\}$
is also strongly mixing with geometric rate, and thus satisfies the
mixing condition $\mcr{A}(a_n)$ described in Remark~\ref{t3.9} with
$r_n=[n^{\delta}]$ for any $0<\delta<1$. Since $\{\varepsilon_k\}$ is a
stationary sequence, we have
 \beqnn
\ar\ar n\mbf{P}\Big(\max_{m\le |k|\le r_n}|\varepsilon_k|> {a}_nx,
|\varepsilon_1|> {a}_nx\Big)\cr
 \ar\ar\qquad
\leq n\sum_{k=m}^{r_n}\Big[\mbf{P}(|\varepsilon_k|> a_nx,
|\varepsilon_1|> a_nx) + \mbf{P}(|\varepsilon_{-k}|> a_nx,
|\varepsilon_1|> a_nx)\Big]\cr
 \ar\ar\qquad
= n\sum_{k=m}^{r_n}\Big[\mbf{P}(|\varepsilon_k|> a_nx, |\varepsilon_1|>
a_nx)+ \mbf{P}(|\varepsilon_{1}|> a_nx, |\varepsilon_{k+2}|>
a_nx)\Big]\cr
 \ar\ar\qquad
\leq 2n\sum_{k=m}^{r_n+2} \mbf{P}(|\varepsilon_k|> a_nx,|\varepsilon_1|>
a_nx).
 \eeqnn
The right hand side tends to zero as $n\to \infty$ by Lemma~\ref{t4.1}.
By Theorem~\ref{t3.16} we have, as $n\to \infty$,
 \beqnn
\mbf{P}\Big(\max_{m\le |k|\le r_n}|\varepsilon_k|> {a}_nx\Big|
|\varepsilon_1|> {a}_nx\Big) \to 0.
 \eeqnn
By Theorem~\ref{t3.16} we have $n\mbf{P}(\|\mbf{I}_1\|>
(cn)^{1/\alpha})\to 1$ as $n\to \infty$, where
 \beqnn
c = \frac{a\sigma^\alpha(1-e^{-\alpha b})} {\alpha^3b^2 \Gamma(-\alpha)}.
 \eeqnn
By Theorem~4.5 in Basrak and Segers (2009), we have (\ref{4.2}) with the
Laplace functional $\mbf{E}[e^{-\eta(f)}]$ given by
 \beqnn
\exp\bigg\{-\frac{1}{\mbf{E}(G)}\int_0^\infty \mbf{E}\Big[\Big(1 -
\exp\Big\{-f\Big(c^{1/\alpha}v,\frac{c^{1/\alpha}v}
{1+X_0}\Big)\Big\}\Big) G\Big]d(-v^{-\alpha})\bigg\}.
 \eeqnn
This clearly coincides with (\ref{4.3}). \qed

Based on the above theorem, we now study the convergence of some partial
sums associated with the sequence $\{\mbf{I}_k\}$ defined by (\ref{3.4}).
To do so, let us introduce some notation. For any $B\in
\mcr{B}(\mbb{R}_+)$ define
 \beqlb\label{4.4}
U_{1,n}(B)=\sum_{k=1}^n\varepsilon_k 1_B(|\varepsilon_k|),
 \quad
U_{2,n}(B) = \sum_{k=1}^n \frac{\varepsilon_k}{1+X_{k-1}}
1_B\Big(\Big|\frac{\varepsilon_k}{1+X_{k-1}}\Big|\Big).
 \eeqlb
Then we define $\tilde{U}_{j,n}(B) = U_{j,n}(B) - \mbf{E}[U_{j,n}(B)]$
for $j=1,2$.

\blemma\label{t4.3} For any $\delta>0$ we have
 \beqnn
\lim_{ z\to0}\limsup_{n\to \infty} \mbf{P}\big(a_n^{-1}
|\tilde{U}_{1,n}(0,a_n z]|> \delta\big)=0.
 \eeqnn
\elemma

\proof Since $\mbf{E}(\varepsilon_k) = \mbf{E}(\varepsilon_k|
\mcr{F}_{k-1}) = 0$, we have
 \beqnn
a_n^{-1}\tilde{U}_{1,n}(0,a_n z]
 \ar=\ar
a_n^{-1}\sum_{k=1}^n [\varepsilon_k1_{\{|\varepsilon_k|\le a_n z\}} -
\mbf{E}(\varepsilon_k1_{\{|\varepsilon_k|\le a_n z\}})]\cr
 \ar=\ar
a_n^{-1}\sum_{k=1}^n [\varepsilon_k1_{\{|\varepsilon_k|\le a_n z\}} -
\mbf{E}(\varepsilon_k1_{\{|\varepsilon_k|\le a_n z\}}|\mcr{F}_{k-1})]\cr
 \ar\ar
-\, a_n^{-1}\sum_{k=1}^n \Big[\mbf{E}(\varepsilon_k 1_{\{|\varepsilon_k|>
a_n z\}}|\mcr{F}_{k-1}) -\mbf{E}(\varepsilon_k1_{\{|\varepsilon_k|> a_n
z\}})\Big].
 \eeqnn
Let $J_1$ and $J_2$ denote the two terms on the right-hand side. By
Theorem~\ref{t3.16} one can see that $\varepsilon_1^2$ is regularly
varying with index $\alpha/2$. Then by Proposition~\ref{t3.2} it follows
that, as $n\to \infty$,
 \beqnn
\mbf{Var}(J_1)
 \ar=\ar
a_n^{-2}\sum_{k=1}^n\mbf{E}\big\{\big[\varepsilon_k1_{\{|\varepsilon_k|\le
a_n z\}} - \mbf{E}(\varepsilon_k1_{\{|\varepsilon_k|\le a_n
z\}}|\mcr{F}_{k-1})\big]^2\big\}\cr
 \ar\le\ar
na_n^{-2}\mbf{E}(\varepsilon_1^21_{\{|\varepsilon_1|\le a_n z\}})
 =
na_n^{-2}\mbf{E}(\varepsilon_1^21_{\{\varepsilon_1^2\le a_n^2 z^2\}})
\sim Cz^{2-\alpha},
 \eeqnn
which goes to zero as $z\to 0$. Now we discuss the asymptotics of $J_2$.
Observe that, for $u>\gamma x+\rho$, we have $|X_1-\gamma x-\rho|> u$ if
and only if $X_1> u+\gamma x+\rho$. It follows that
 \beqnn
~ u^{\alpha-1}\mbf{E}_x\big(|X_1-\varepsilon_1|1_{\{|\varepsilon_1|>
u\}}\big)
 \ar=\ar
u^{\alpha-1}\mbf{E}_x\big(|\gamma x+\rho|1_{\{|X_1-\gamma x-\rho|>
u\}}\big) \cr
 \ar=\ar
(\gamma x+\rho)u^{\alpha-1}\mbf{P}_x(X_1>u+\gamma x+\rho).
 \eeqnn
Using Proposition~\ref{t3.10} we see the right-hand side tends to zero
uniformly in $x\in [0,K]$ as $u\to \infty$. By Proposition~\ref{t3.11},
we have
 \beqnn
\lim_{u\to \infty}u^{\alpha-1}\mbf{E}_x\big[|\varepsilon_1|
1_{\{|\varepsilon_1|> u\}}\big]
 \ar=\ar
\lim_{u\to \infty}u^{\alpha-1}\mbf{E}_x\big[\varepsilon_1
1_{\{|X_1-\gamma x-\rho|> u\}}\big] \cr
 \ar=\ar
\lim_{u\to \infty} u^{\alpha-1}\mbf{E}_x[X_11_{\{X_1> u+\gamma x+\rho\}}]
\cr
 \ar=\ar
\frac{\sigma^\alpha}{(\alpha-1)\Gamma(-\alpha)} (q_\alpha+p_\alpha x),
 \eeqnn
and the convergence is uniform in $x\in [0,K]$. It follows that, as $n\to
\infty$, we have almost surely
 \beqlb\label{4.5}
\ar\ar a_n^{-1}\sum_{k=1}^n1_{\{X_{k-1}\le K\}}
\mbf{E}\big[\varepsilon_{k} 1_{\{|\varepsilon_k|>
a_nz\}}\big|\mcr{F}_{k-1} \big] \cr
 \ar\ar\qquad
= a_n^{-1}\sum_{k=1}^n 1_{\{X_{k-1}\le K\}} \mbf{E}_{X_{k-1}}
\big[\varepsilon_1 1_{\{|\varepsilon_1|> a_n z\}}\big] \cr
 \ar\ar\qquad
= \frac{\sigma^\alpha z^{1-\alpha}}{(\alpha-1)\Gamma(-\alpha)n}
\sum_{k=1}^n 1_{\{X_{k-1}\le K\}} (q_\alpha+p_\alpha X_{k-1}) + o(1) \cr
 \ar\ar\qquad
= \frac{\sigma^\alpha z^{1-\alpha}}{(\alpha-1)\Gamma(-\alpha)n}
\sum_{k=1}^n 1_{\{X_{k-1}\le K\}}(q_\alpha+p_\alpha X_{k-1}) + o(1) \cr
 \ar\ar\qquad
= \frac{\sigma^\alpha z^{1-\alpha}}{(\alpha-1)\Gamma(-\alpha)}
\mbf{E}\big[1_{\{X_{0}\le K\}}(q_\alpha+p_\alpha X_0)\big] + o(1),
 \eeqlb
where the last equality holds the ergodic theorem. Similarly, we have
 \beqlb\label{4.6}
\ar\ar na_n^{-1}\mbf{E}\big[1_{\{X_{0}\le K\}}\varepsilon_1
1_{\{|\varepsilon_1|> a_n z\}}\big] \cr
 \ar\ar\qquad
= \frac{\sigma^\alpha z^{1-\alpha}}{(\alpha-1)\Gamma(-\alpha)}
\mbf{E}\big[1_{\{X_{0}\le K\}}(q_\alpha+p_\alpha X_0)\big] + o(1).
 \eeqlb
Then (\ref{4.5}) and (\ref{4.6}) cancel asymptotically as $n\to \infty$.
Observe that
 \beqnn
\mbf{P}(1_{\{X_0>K\}}|\varepsilon_1|> u)
 \ar=\ar
\mbf{P}\Big(\sigma\int_0^11_{\{X_0>K\}}e^{-b(1-s)} X_{s-}^{1/\alpha}
dZ_s> u\Big) \cr
 \ar\ar
+\,\mbf{P}\Big(\sigma\int_0^11_{\{X_0>K\}}e^{-b(1-s)} X_{s-}^{1/\alpha}
d(-Z_s)> u\Big).
 \eeqnn
By Remark~\ref{t3.8}, as $u\to \infty$,
 \beqnn
\mbf{P}(1_{\{X_0>K\}}|\varepsilon_1|> u)
 \sim
C(K)\big[\mbf{P}(Z_1>u) + \mbf{P}(-Z_1>u)\big]
 =
C(K)u^{-\alpha},
 \eeqnn
where
 \beqnn
C(K) = \sigma^{\alpha}\mbf{E}\Big[\int_0^11_{\{X_0>K\}} e^{-\alpha
b(1-s)} X_sds\Big]
 \le
C_1\mbf{E}[1_{\{X_0>K\}}(1+X_0)].
 \eeqnn
Then by Proposition~\ref{t3.2}, as $n\to \infty$,
 \beqnn
\ar\ar a_n^{-1}\mbf{E}\Big\{\sum_{k=1}^n \big[1_{\{X_{k-1}>K\}}
\mbf{E}(|\varepsilon_k| 1_{\{|\varepsilon_k|> a_n z\}} |\mcr{F}_{k-1})
\cr
 \ar\ar\qquad\qquad\qquad\qquad
+\, \mbf{E}(1_{\{X_{k-1}> K\}} |\varepsilon_k| 1_{\{|\varepsilon_k|>
a_nz\}})\big]\Big\} \cr
 \ar\ar\qquad
=\, 2a_n^{-1}\sum_{k=1}^n \mbf{E}\big(1_{\{X_{k-1}> K\}} |\varepsilon_k|
1_{\{1_{\{X_{k-1}> K\}}|\varepsilon_k|> a_nz\}}\big) \ccr
 \ar\ar\qquad
= a_n^{-1}n\mbf{E}\big[1_{\{X_0> K\}}|\varepsilon_1| 1_{\{1_{\{X_0>
K\}}|\varepsilon_1|> a_nz\}}\big] \ccr
 \ar\ar\qquad
= C_2n\mbf{P}\big(1_{\{X_0> K\}}|\varepsilon_1|> a_nz\big)
 =
C_2C(K) z^{-\alpha}.
 \eeqnn
The right hand side goes to zero as $K\to \infty$. That gives the desired
result. \qed

\blemma\label{t4.4} For any $\delta>0$ we have
 \beqnn
\lim_{z\to 0}\limsup_{n\to \infty} \mbf{P}\big(a_n^{-1}
|\tilde{U}_{2,n}(0,a_n z]|> \delta\big)=0.
 \eeqnn
\elemma

\proof It is simple to see that
 \beqnn
u^{\alpha-1}\mbf{E}_x\Big[\frac{\varepsilon_1}{x+1}
1_{\{|\varepsilon_1|>(x+1)u\}}\Big]
 =
\frac{[u(x+1)]^{\alpha-1}}{(x+1)^{\alpha}} \mbf{E}_x[\varepsilon_1
1_{\{|\varepsilon_1|>(x+1)u\}}],
 \eeqnn
where $u(x+1)>u$ and $(x+1)^{-\alpha}\le 1$. Thus as $u\to \infty$,
uniformly for $x\in[0,K]$,
 \beqnn
u^{\alpha-1}\mbf{E}_x\Big[\frac{\varepsilon_1}{x+1}1_{\{|\varepsilon_1|>
(x+1)u\}}\Big]\to \frac{\sigma^\alpha(q_\alpha+p_\alpha
x)}{(\alpha-1)\Gamma(-\alpha)(x+1)^{\alpha}}.
 \eeqnn
The remaining argument is similar to the proof of Lemma~\ref{t4.3}. \qed

\ttheorem\label{t4.5} Let $U_{1,n} = U_{1,n}(0,\infty)$ and $U_{2,n} =
U_{2,n}(0,\infty)$. Then we have, as $n\to \infty$,
 \beqnn
a_n^{-1}(U_{1,n},U_{2,n})\overset{\rm d}{\longrightarrow} (U_1,U_2)
\quad\mbox{on}\quad \mbb{R}^2,
 \eeqnn
where $(U_1,U_2)$ is the $\alpha$-stable random vector with
characteristic function given by
 \beqlb\label{4.7}
\ar\ar\mbf{E}\big[\exp\{i(\lambda_1U_1 + \lambda_2U_2)\}\big] \ccr
 \ar\ar\qquad
= \exp\Big\{\frac{\sigma^\alpha}{\alpha} \mbf{E}\Big[\Big(\lambda_1 +
\frac{\lambda_2} {1+X_0}\Big)^\alpha \big(q_\alpha + p_\alpha
X_0\big)\Big]e^{-i\pi\alpha/2}\Big\},
 \eeqlb
and $p_\alpha = p_\alpha(1)$ and $q_\alpha = q_\alpha(1)$ are defined by
(\ref{3.2}). \etheorem

\proof Fix $z>0$ and $\lambda=(\lambda_1,\lambda_2)\in\mbb{R}^2$, and
define the function on $\mbb{R}^2$ by $f_{\lambda,z}(x_1,x_2)=
\lambda_1x_1 1_{\{|x_1|> z\}} + \lambda_2x_2 1_{\{|x_2|> z\}}$. Then we
have
 \beqnn
\eta_n(f_{\lambda,z})
 =
a_n^{-1}\sum_{j=1}^2 \lambda_j U_{j,n}(a_nz,\infty).
 \eeqnn
It is easy to see that the mapping from $M(\mbb{R}^2)$ into $\mbb{R}$
defined by
 \beqnn
N:= \sum_{k=1}^\infty\delta_{(x_{1,k},x_{2,k})} \mapsto N(f_{\lambda,z})
 \eeqnn
is a.s.\ continuous with respect to the distribution of the limit point
process $\eta$ in Theorem~\ref{t4.2}. By the continuous mapping theorem,
as $n\to \infty$, , we have $\eta_n(f_{\lambda,z})\overset{\rm d}
{\longrightarrow} \eta(f_{\lambda,z})$, and hence
 \beqnn
\mbf{E}\Big[\exp\Big\{ia_n^{-1}\sum_{j=1}^2 \lambda_j
U_{j,n}(a_nz,\infty)\Big\}\Big]
 =
\mbf{E}\big[\exp\{i\eta_n(f_{\lambda,z})\}\big]
 \to
\mbf{E}\big[\exp(i\eta(f_{\lambda,z}))\big],
 \eeqnn
where the right-hand side is given by
 \beqnn
\exp\bigg\{\frac{1}{\alpha\Gamma(-\alpha)}\int_0^\infty
\mbf{E}\Big[\Big(\exp\Big\{{i\lambda_1y1_{\{y>z\}} + \frac{i\lambda_2y
1_{\{y>z(1+X_0)\}}} {1+X_0}}\Big\} - 1\Big)G\Big]
\frac{dy}{y^{\alpha+1}}\bigg\}.
 \eeqnn
By Proposition~\ref{t3.2} and Theorem~\ref{t3.16}, as $n\to \infty$,
 \beqnn
a_n^{-1}\mbf{E}[U_{1,n}(a_n z,\infty)]
 \ar=\ar
na_n^{-1} \mbf{E}(\varepsilon_1 1_{\{|\varepsilon_1|> a_n z\}})
 \sim
\frac{\alpha nz}{\alpha-1}\mbf{P}(|\varepsilon_1|> a_n z) \cr
 \ar\sim\ar
\frac{1}{\alpha(\alpha-1)\Gamma(-\alpha)}\mbf{E}(G)z^{1-\alpha} \cr
 \ar=\ar
\frac{1}{\alpha\Gamma(-\alpha)}\int_0^\infty \mbf{E}(G)y1_{\{y>z\}}
\frac{dy}{y^{\alpha+1}}.
 \eeqnn
By Proposition~\ref{t3.2} and Remarks~\ref{t3.7} and~\ref{t3.8}, as $n\to
\infty$,
 \beqnn
a_n^{-1}\mbf{E}[U_{2,n}(a_n z,\infty)]
 \ar=\ar
na_n^{-1}\mbf{E}\Big[\frac{\varepsilon_1} {1+X_0}
1_{\big\{\big|\frac{\varepsilon_1}{{1+X_0}}\big|> a_n z\big\}}\Big]\cr
 \ar\sim\ar
\frac{\alpha nz}{\alpha-1}
\mbf{P}\Big(\Big|\frac{\varepsilon_1}{{1+X_0}}\Big|> a_n z\Big)\cr
 \ar\sim\ar
\frac{1}{\alpha(\alpha-1)\Gamma(-\alpha)}
\mbf{E}\Big[\frac{G}{(1+X_0)^\alpha}\Big] z^{1-\alpha} \cr
 \ar=\ar
\frac{1}{\alpha\Gamma(-\alpha)}\int_0^\infty
\mbf{E}\Big[\frac{G}{(1+X_0)^\alpha}\Big] y1_{\{y>z\}}
\frac{dy}{y^{\alpha+1}}.
 \eeqnn
Consequently, as $n\to \infty$,
 \beqnn
\mbf{E}\Big[\exp\Big\{ia_n^{-1}\sum_{j=1}^2 \lambda_j
\tilde{U}_{j,n}(a_nz,\infty)\Big\}\Big]
 =
\mbf{E}\big[\exp\{i\tilde{\eta}_n(f_{\lambda,z})\}\big]
 \eeqnn
converges to
 \beqnn
\ar\ar\exp\bigg\{\frac{1}{\alpha\Gamma(-\alpha)}\int_0^\infty
\mbf{E}\Big[\Big(\exp\Big\{i\lambda_1y1_{\{y>z\}} + \frac{i\lambda_2y}
{1+X_0}1_{\{y>z(1+X_0)\}}\Big\} - 1 \cr
 \ar\ar\qquad\qquad\qquad\qquad\qquad
-\, i\lambda_1y1_{\{y>z\}} - \frac{i\lambda_2y} {1+X_0}1_{\{y>
z(1+X_0)\}}\Big)G\Big] \frac{dy}{y^{\alpha+1}}\bigg\}.
 \eeqnn
As $z\to 0$, the above quality tends to
 \beqnn
\ar\ar \exp\bigg\{\frac{1}{\alpha\Gamma(-\alpha)}\int_0^\infty
\mbf{E}\Big[\Big(e^{i\lambda_1y+\frac{i\lambda_2y}{1+X_0}}-1-
i\lambda_1y-\frac{i\lambda_2y}{1+X_0}\Big)G\Big]
\frac{dy}{y^{\alpha+1}}\bigg\} \cr
 \ar\ar\qquad
= \exp\bigg\{\frac{1}{\alpha\Gamma(-\alpha)}\int_0^\infty \big(e^{iz}-1-
iz\big) \mbf{E}\Big[\Big(\lambda_1 + \frac{\lambda_2} {1+X_0}\Big)^\alpha
G\Big] \frac{dz}{z^{\alpha+1}}\bigg\}.
 \eeqnn
By Corollary~14.11 of Sato (1999) and (\ref{3.3}) one can see this
coincides with (\ref{4.7}). Since $\mbf{E}(U_{j,n}) =
\mbf{E}[U_{j,n}(0,\infty)] = 0$, by the above calculations and
Lemmas~\ref{t4.3} and~\ref{t4.4}, as $n\to \infty$,
 \beqnn
\mbf{E}\Big[\exp\Big\{ia_n^{-1}\sum_{j=1}^2 \lambda_j U_{j,n}\Big\}\Big]
 =
\mbf{E}\Big[\exp\Big\{ia_n^{-1}\sum_{j=1}^2 \lambda_j
\tilde{U}_{j,n}(0,\infty)\Big\}\Big]
 \eeqnn
converges to (\ref{4.7}). That gives the desired result. \qed

\blemma\label{t4.6} Let $r_n=[n^\delta]$ with $0<\delta<1$. Then we have
 \beqnn
\lim_{m\to \infty} \limsup_{n\to \infty} n\mbf{P}\Big(\max_{-r_n\leq
k\leq -m}\|\mbf{H}_k\|>c_nx, X_0> a_nx\Big) = 0.
 \eeqnn
\elemma

\proof Since $\{(\mbf{H}_k,X_k): k\in \mbb{Z}\}$ is a stationary
sequence, by (\ref{1.5}), it is easy to see
 \beqnn
\ar\ar n\mbf{P}\Big(\max_{-r_n\leq k\leq -m}\|\mbf{H}_k\|> c_nx, X_0>
a_nx\Big) \cr
 \ar\ar\qquad
= n\mbf{P}\Big(\max_{m-r_n\leq k\leq 0}\|\mbf{H}_k\|>c_nx, X_m> a_nx\Big)
\cr
 \ar\ar\qquad
\leq n\mbf{P}\Big(e^{-bm}X_0 + a\int_0^m e^{-b(m-s)} ds> a_nx/2\Big) \cr
 \ar\ar\qquad\quad
+\, n\mbf{P}\Big(A_n, \sigma\Big|\int_0^m e^{-b(m-s)} X_{s-}^{1/\alpha}
dZ_s\Big|> a_nx/2\Big),
 \eeqnn
where
 \beqnn
A_n = \Big(\max_{-r_n\le k\leq 0}\|\mbf{H}_k\|> c_nx\Big).
 \eeqnn
By Theorem~\ref{t3.19} it is easy to see that
 \beqnn
\lim_{n\rightarrow\infty}\mbf{P}(A_n)
 \ar=\ar
\lim_{n\rightarrow\infty}\sum_{k=0}^{r_n} \mbf{P}(\|\mbf{H}_{-k}\|> c_nx)
\cr
 \ar=\ar
\lim_{n\rightarrow\infty}(r_n+1) \mbf{P}(\|\mbf{H}_0\|> c_nx) = 0.
 \eeqnn
Then Lemma~\ref{t3.14} implies that
 \beqnn
\lim_{n\rightarrow\infty} n\mbf{P}\Big(A_n, \sigma\Big|\int_0^m
e^{-b(m-s)} X_{s-}^{1/\alpha} dZ_s\Big|> a_nx/2\Big) = 0.
 \eeqnn
From Proposition~\ref{t3.12} it follows that
 \beqnn
\lim_{n\rightarrow\infty} n\mbf{P}\Big(e^{-bm}X_0 + a\int_0^m e^{-b(m-s)}
ds> a_nx/2\Big)
 =
Ce^{-\alpha bm}x^{-\alpha},
 \eeqnn
which goes to zero as $m\rightarrow\infty$. Then we have the desired
result. \qed

\blemma\label{t4.7} There exists $\delta\in (0,1)$ so that for
$r_n=[n^\delta]$ we have
 \beqnn
\lim_{m\to \infty}\limsup_{n\to \infty} n\mbf{P}\Big(\max_{m\leq k\leq
r_n}\|\mbf{H}_k\|>c_nx, X_0> a_nx\Big) = 0.
 \eeqnn
\elemma

\proof Recall that $X_0$ is regularly varying with index $\alpha$. It is
easy to see that
 \beqnn
\ar\ar n\mbf{P}\Big(\max_{m\leq k\leq r_n}\|\mbf{H}_k\|>c_nx, X_0>
a_nx\Big) \cr
 \ar\ar\qquad
\le n\sum_{k=m}^{r_n}\mbf{P}(\|\mbf{H}_k\|>c_nx, X_0> a_nx)\cr
 \ar\ar\qquad
\le n\sum_{k=m}^{r_n}\mbf{P}(\|\bar{\mbf{H}}_k-\mbf{H}_k\|> c_nx/2) \cr
 \ar\ar\qquad\quad
+\, n\sum_{k=m}^{r_n}\mbf{P}(\|\bar{\mbf{H}}_k\|> c_nx/2, X_0> a_nx).
 \eeqnn
Let $J_1$ and $J_2$ denote the two terms on the right-hand side. We can
choose $r_n=[n^\delta]$ for sufficiently small $\delta\in (0,1)$, and
thus
 \beqnn
\limsup_{n\to \infty}J_1
 \ar\le\ar
\limsup_{n\to \infty} \frac{Cnr_n}{x^rc_n^r}
\mbf{E}(\|\bar{\mbf{H}}_1-\mbf{H}_1\|^r) \cr
 \ar\le\ar
\limsup_{n\to \infty} \frac{Cn^{1+\delta}}{x^rc_n^{r}}
\mbf{E}(\|\bar{\mbf{H}}_1-\mbf{H}_1\|^r)
 = 0.
 \eeqnn
By Proposition~\ref{t3.2}, we have
 \beqnn
\mbf{E}[X_01_{\{X_0> a_nx\}}]
 \sim
\frac{\alpha a_nx}{\alpha-1} \mbf{P}(X_0> a_nx)
 \sim
C(a_nx)^{1-\alpha}.
 \eeqnn
By Remark~\ref{t3.17}, we have $E[|1\vee
V_k|^{\alpha/(\alpha+1)}]<\infty$. Note that $(X_0, X_{k-1})$ is
independent of $V_k$ for $k\ge 2$. Then for some constant
$\delta\in(0,1/\alpha)$,
 \beqnn
J_2\ar\le\ar\frac{2^{\alpha/(\alpha+1)}n}{a_nx^{\alpha/(\alpha+1)}}
\sum_{k=m}^{r_n} \mbf{E}\big[X_{k-1}|1\vee V_k|^{\alpha/(\alpha+1)}; X_0>
a_nx\big]\cr
 \ar\le\ar
\frac{Cn}{a_nx^{\alpha/(\alpha+1)}} \sum_{k=m}^{r_n} \mbf{E}\big[|1\vee
V_k|^{\alpha/(\alpha+1)}\big] \mbf{E}\big[1_{\{X_0> a_nx\}}
\mbf{E}_{X_0}(X_{k-1})\big] \cr
 \ar\le\ar
\frac{Cn}{a_nx^{\alpha/(\alpha+1)}}\sum_{k=m}^{r_n}\mbf{E}\Big\{1_{\{X_0>
a_nx\}} \big[X_0e^{-b(k-1)} + ab^{-1}(1-e^{-b(k-1)})\big]\Big\} \cr
 \ar\le\ar
\frac{Cn}{a_nx^{\alpha/(\alpha+1)}}\mbf{E}[X_01_{\{X_0> a_nx\}}]
\sum_{k=m}^{r_n} e^{-b(k-1)} + \frac{Cnr_n}{a_nx^{\alpha/(\alpha+1)}}
\mbf{P}(X_0> a_nx).
 \eeqnn
It follows that
 \beqnn
\limsup_{n\to \infty}J_2\le Cz^{1-\alpha}\sum_{k=m}^{\infty} e^{-b(k-1)},
 \eeqnn
which goes to zero as $m\to \infty$. \qed

\blemma\label{t4.8} Let $\bar{\mbf{H}}_k = X_{k-1}^{(\alpha+1)/\alpha}
(1,V_k)$. Then there exists $\delta\in (0,1)$ so that for
$r_n=[n^\delta]$ we have
 \beqnn
\lim_{m\to \infty}\limsup_{n\to \infty} n\mbf{P}\Big(\max_{m\le |k|\le
r_n}\|\mbf{H}_k\|> c_nx, \|\bar{\mbf{H}}_1\|>c_nx\Big) = 0.
 \eeqnn
\elemma

\proof For any $K>1$, we have
 \beqnn
\ar\ar n\mbf{P}\Big(\max_{m\le |k|\le r_n}\|\mbf{H}_k\|> c_nx,
\|\bar{\mbf{H}}_1\|>c_nx\Big) \cr
 \ar\ar\qquad
\leq n\mbf{P}\big(X_0^{(\alpha+1)/\alpha}|V_1|1_{\{|V_1|>K\}}> c_nx/2\big)
\cr
 \ar\ar\qquad\quad
+\, n\mbf{P}\Big(\max_{m\leq k\leq r_n}\|\mbf{H}_k\|>c_nx,
KX_0^{(\alpha+1)/\alpha}>c_nx/2\Big) \cr
 \ar\ar\qquad\quad
+\,n\mbf{P}\Big(\max_{-r_n\leq k\leq -m}\|\mbf{H}_k\|>c_nx,
KX_0^{(\alpha+1)/\alpha}>c_nx/2\Big).
 \eeqnn
Observe that $X_0^{(\alpha+1)/\alpha}$ is regularly varying with index
$\alpha^2/(\alpha+1)$. By Remark~\ref{t3.17}, we have $\mbf{E}[|V_1|^b]<
\infty$ for some $b>\alpha^2/(\alpha+1)$. It follows from Breiman's Lemma
that
 \beqnn
\ar\ar\lim_{n\to \infty} n\mbf{P}\big(X_0^{(\alpha+1)/\alpha}|V_1|
1_{\{|V_1|>K\}}> c_nx/2\big) \cr
 \ar\ar\qquad
= \lim_{n\to \infty} n\mbf{E}(|V_1|^{\alpha^2/(\alpha+1)} 1_{\{|V_1|>
K\}}) \mbf{P}(X_0^{(\alpha+1)/\alpha}> c_nx/2) \cr
 \ar\ar\qquad
= C x^{-\alpha^2/(\alpha+1)}\mbf{E}(|V_1|^{\alpha^2/(\alpha+1)}
1_{\{|V_1|> K\}}).
 \eeqnn
The right-hand side goes to zero as $K\to \infty$. Then the result
follows by Lemmas~\ref{t4.6} and~\ref{t4.7}. \qed

\ttheorem\label{t4.9} Let $\{V_j\}$ be defined by (\ref{3.14}). Then we
have, as $n\to \infty$,
 \beqlb\label{4.8}
\xi_n := \sum_{k=1}^n\delta_{c_n^{-1}\mbf{H}_k} \overset{\rm
d}{\longrightarrow} \xi \quad\mbox{on}\quad M(\bar{\mbb{R}}_0^2),
 \eeqlb
where $\xi$ is a point process on $\bar{\mbb{R}}_0^2$ with Laplace
functional $\mbf{E}[e^{-\xi(f)}]$, $f\in C_0^+(\bar{\mbb{R}}_0^2)$ given
by
 \beqlb\label{4.9}
\ar\ar \exp\bigg\{-\frac{a\sigma^\alpha}{\alpha^2b^2\Gamma(-\alpha)}
\int_0^\infty \mbf{E}\Big(1-\exp\big\{-f\big(y^{(\alpha+1)/\alpha}
(1,V_1)\big)\big\}\Big) \cr
 \ar\ar\qquad\quad
\mbf{E}\Big[\exp\Big\{-\sum_{j=2}^\infty f\big(y^{(\alpha+1)/\alpha}
e^{-b(j-1)(\alpha+1)/\alpha} (1,V_j)\big)\Big\}\Big]
\frac{dy}{y^{\alpha+1}}\bigg\}.
 \eeqlb
 \etheorem

\proof This proof is based on Theorem~4.5 of Basrak and Segers (2009)
similarly as the proof of Theorem~\ref{t4.2}. It is easy to see that
 \beqnn
\ar\ar n\mbf{P}\Big(\max_{m\leq |k|\leq r_n}\|\mbf{H}_k\|>
c_nx,\|\mbf{H}_1\|> c_nx\Big)\cr
 \ar\ar\qquad
\leq n\mbf{P}\big(\|\bar{\mbf{H}}_1-\mbf{H}_1\|>c_nx/2\big) \cr
 \ar\ar\qquad\quad
+\, n\mbf{P}\Big(\max_{m\leq |k|\leq r_n}\|\mbf{H}_k\|>c_nx,
\|\bar{\mbf{H}}_1\|> c_nx/2\Big).
 \eeqnn
By Lemma~\ref{t3.18}, we have $\mbf{E}[\|\mbf{H}_k -
\bar{\mbf{H}}_k\|^r]< \infty$ for some $r>\alpha^2/(\alpha+1)$. Then
Markov's inequality implies that
 \beqnn
\limsup_{n\to \infty} n\mbf{P}\big(\|\bar{\mbf{H}}_1-\mbf{H}_1\|>
c_nx/2\big)
 \le
\limsup_{n\to \infty} nc_n^{-r}\mbf{E}\big[\|\mbf{H}_1-
\bar{\mbf{H}}_1\|^r\big]=0.
 \eeqnn
By Lemma~\ref{t4.8} we have
 \beqnn
\lim_{m\to \infty}\limsup_{n\to \infty}n \mbf{P}\Big(\max_{m\le |k|\le
r_n}\|\mbf{H}_k\|>c_nx,\|\mbf{H}_1\|>c_nx\Big)=0,
 \eeqnn
where $r_n=[n^\delta]$ for some $\delta\in (0,1)$. Let
 \beqnn
h = \frac{a\sigma^\alpha}{\alpha^3b^2\Gamma(-\alpha)}\mbf{E}[1\vee
|V_1|^{\alpha^2/(\alpha+1)}].
 \eeqnn
By Theorem~\ref{t3.19} we have, as $n\to \infty$,
 \beqnn
n\mbf{P}\{\|\mbf{H}_1\|>(hn)^{(\alpha+1)/\alpha^2}\}
 =
n\mbf{P}\{\|\mbf{H}_1\|>h^{(\alpha+1)/\alpha^2}c_n\}\to 1.
 \eeqnn
Observe also that
 \beqnn
\xi_n(f) = \sum_{k=1}^nf(c_n^{-1}\mbf{H}_k)
 =
\sum_{k=1}^nf(h^{(\alpha+1)/\alpha^2}(h^{(\alpha+1)/\alpha^2}
c_n)^{-1}\mbf{H}_k).
 \eeqnn
Let $\mbf{\Theta}_i$ be defined as in Theorem~\ref{t3.19}. Then we can
use Theorem~4.5 of Basrak and Segers (2009) to obtain (\ref{4.8}) with
$\mbf{E}[e^{-\xi(f)}]$ given by
 \beqnn
\ar\ar \exp\bigg\{-\frac{1}{\mbf{E}(1\vee |V_1|^{\alpha^2/(\alpha+1)})}
\int_0^\infty \mbf{E}\Big[ \exp\Big\{-\sum_{j=2}^\infty
f(h^{(\alpha+1)/\alpha^2} v\mbf{\Theta}_i)\Big\}\cr
 \ar\ar\qquad\quad
\Big(1-\exp\big\{-f(h^{(\alpha+1)/\alpha^2} v \mbf{\Theta}_1)\big\}\Big)
(1\vee |V_1|^{\alpha^2/(\alpha+1)})\Big]
d(-v^{-\alpha^2/(\alpha+1)})\bigg\} \cr
 \ar\ar\qquad
= \exp\bigg\{-\frac{h}{\mbf{E}(1\vee |V_1|^{\alpha^2/ (\alpha+1)})}
\int_0^\infty \mbf{E}\Big[ \exp\Big\{-\sum_{j=2}^\infty
f(u^{(\alpha+1)/\alpha} \mbf{\Theta}_i)\Big\}\cr
 \ar\ar\qquad\quad
\Big(1-\exp\big\{- f(u^{(\alpha+1)/\alpha} \mbf{\Theta}_1)\big\}\Big)
(1\vee |V_1|^{\alpha^2/(\alpha+1)})\Big] d(-u^{-\alpha})\bigg\} \cr
 \ar\ar\qquad
= \exp\bigg\{-\frac{a\sigma^\alpha}{\alpha^3b^2\Gamma(-\alpha)}
\int_0^\infty \mbf{E}\Big[ \exp\Big\{-\sum_{j=2}^\infty
f(u^{(\alpha+1)/\alpha} \mbf{\Theta}_i)\Big\}\cr
 \ar\ar\qquad\quad
\Big(1-\exp\big\{- f(u^{(\alpha+1)/\alpha} \mbf{\Theta}_1)\big\}\Big)
(1\vee |V_1|^{\alpha^2/(\alpha+1)})\Big] d(-u^{-\alpha})\bigg\} \cr
 \ar\ar\qquad
= \exp\bigg\{-\frac{a\sigma^\alpha}{\alpha^3b^2\Gamma(-\alpha)}
\int_0^\infty \mbf{E}\Big[\exp\Big\{-\sum_{j=2}^\infty
f(y^{(\alpha+1)/\alpha}(1\vee |V_1|)\mbf{\Theta}_i)\Big\}\cr
 \ar\ar\qquad\quad
\Big(1-\exp\big\{-f(y^{(\alpha+1)/\alpha}(1\vee |V_1|)
\mbf{\Theta}_1)\big\}\Big)\Big]d(-y^{-\alpha})\bigg\},
 \eeqnn
which can be rewritten as (\ref{4.9}). \qed

From the above theorem, we can derive some limit theorem of partial sums
associated with the sequence $\{\mbf{H}_k\}$ defined by (\ref{3.4}). For
$B\in \mcr{B}(\mbb{R}_+)$ define
 \beqlb\label{4.10}
S_{1,n}(B)=\sum_{k=1}^nX^2_{k-1} 1_B(X_{k-1}),
 \quad
S_{2,n}(B) = \sum_{k=1}^n X_{k-1}\varepsilon_k
1_B(|X_{k-1}\varepsilon_k|).
 \eeqlb

\blemma\label{t4.10} For any $\delta> 0$ we have
 \beqnn
\lim_{z\to 0}\limsup_{n\to \infty} \mbf{P}\big(c_n^{-2} |S_{1,n}(0,c_n
z)|> \delta\big)=0.
 \eeqnn
\elemma

\proof By Theorem~\ref{t3.15}, it is easy to see that $X_0^2$ is
regularly varying with index $\alpha/2< 1$. Using Proposition~\ref{t3.2}
and Theorem~\ref{t3.15}, we have, as $n\rightarrow\infty$,
 \beqnn
\mbf{E}\big[c_n^{-2}S_{1,n}(0,c_nz)\big]
 =
\frac{1}{c_n^2}\sum_{k=1}^n\mbf{E}\big[X_{k-1}^2 1_{\{X_{k-1}<
c_nz\}}\big]
 \sim
\frac{n\alpha z^2}{2-\alpha}\mbf{P}(X_0> c_nz)
 \sim
Cz^{2-\alpha}.
 \eeqnn
The right-hand side tends to zero as $z\rightarrow 0$. Then we have the
desired result. \qed

\blemma\label{t4.11} Suppose that $1<\alpha< (1+\sqrt{5})/2$. Then for
any $\delta> 0$ we have
 \beqnn
\lim_{z\to 0}\limsup_{n\to \infty} \mbf{P}\big(c_n^{-1} |S_{2,n}(0,c_n
z)|> \delta\big)=0.
 \eeqnn
\elemma

\proof By Theorem~\ref{t3.19}, we see $X_0\varepsilon_1$ is regularly
varying with index $\alpha^2/(\alpha+1)$. Under the condition $1<\alpha<
(1+\sqrt{5})/2$, we have $\alpha^2/(\alpha+1)< 1$. By
Proposition~\ref{t3.2} and Theorem~\ref{t3.19}, as $n\rightarrow\infty$,
 \beqnn
\mbf{E}\big[c_n^{-1}|S_{2,n}(0,c_n z)|\big]
 \ar\leq\ar
\frac{1}{c_n}\sum_{k=1}^n\mbf{E}\big[|X_{k-1} \varepsilon_k|
1_{\{|X_{k-1} \varepsilon_k|< c_nz\}}\big]
 =
\frac{n}{c_n}\mbf{E}\big[|X_0\varepsilon_1| 1_{\{|X_0\varepsilon_1|<
c_nz\}}\big] \cr
 \ar\sim\ar
\frac{\alpha^2nz}{\alpha+1-\alpha^2}\mbf{P}(|X_0 \varepsilon_1|> c_nz)
 \sim
Cz^{1-\alpha^2/(\alpha+1)}.
 \eeqnn
The right-hand side tends to zero as $z\rightarrow 0$. That gives the
result; see also Davis and Hsing (1995, p.896). \qed

\ttheorem\label{t4.12} Let $V_1$ be defined by (\ref{3.14}). Let $S_{1,n}
= S_{1,n}(0,\infty)$ and $S_{2,n} = S_{2,n}(0,\infty)$. If $1<\alpha<
(1+\sqrt{5})/2$, then we have, as $n\to \infty$,
 \beqnn
(a_n^{-2}S_{1,n},c_n^{-1}S_{2,n}) \overset{\rm d}{\longrightarrow}
(S_1,S_2) \quad\mbox{on}\quad \mbb{R}^2,
 \eeqnn
where $(S_1,S_2)$ has characteristic function $\mbf{E}[\exp\{i\lambda_1
S_1 + i\lambda_2S_2\}]$ given by
 \beqlb\label{4.11}
\ar\ar~ \exp\bigg\{-\frac{a\sigma^\alpha}{\alpha^2b^2\Gamma(-\alpha)}
\int_0^\infty \mbf{E}\Big(1-\exp\big\{i\lambda_1y^2+i\lambda_2
y^{(\alpha+1)/\alpha} V_1\big\}\Big) \cr
 \ar\ar\qquad\qquad
\mbf{E}\Big[\exp\Big\{\frac{ie^{-2b}\lambda_1 y^2}{1-e^{-2b}} +
\frac{ie^{-b(\alpha+1)/\alpha} \lambda_2y^{(\alpha+1)/\alpha}V_2}
{(1-e^{-b(\alpha+1)})^{1/\alpha}}\Big\}\Big]
\frac{dy}{y^{\alpha+1}}\bigg\}.
 \eeqlb
\etheorem

\proof We first remark that the integral on the right-hand side of
(\ref{4.11}) is well-defined. In fact, by Remarks~\ref{t3.7}
and~\ref{t3.17}, we have, as $x\to \infty$,
 \beqnn
\mbf{P}(V_1\ge x)\sim C_1x^{-\alpha} + o(x^{-\alpha}),
 \qquad
\mbf{P}(V_1\le -x) = o(x^{-\alpha}).
 \eeqnn
By Theorems~8.1.10 and~8.1.11 in Bingham et al.\ (1987), we have, as
$\lambda\to 0$,
 \beqnn
\mbf{E}[1-\cos(\lambda V_1)]\sim C_2\lambda^\alpha + o(\lambda^\alpha),
 \quad
\mbf{E}[\sin(\lambda V_1)]\sim C_3\lambda^\alpha + o(\lambda^\alpha).
 \eeqnn
It follows that $\mbf{E}(1-e^{i\lambda V_1})\sim c\lambda^\alpha$ as
$\lambda\to 0$. Then the integral in (\ref{4.11}) converges. Fix any
$\lambda = (\lambda_1,\lambda_2)\in \mbb{R}^2$ and $z>0$, define the
function on $\mbb{R}_+\times\mbb{R}$ by
 \beqnn
g_{\lambda,z}(x_1,x_2)
 =
\lambda_1x_1^{2\alpha/(\alpha+1)} 1_{\{x_1> z^{(\alpha+1)/\alpha}\}} +
\lambda_2x_21_{\{|x_2|>z\}}.
 \eeqnn
It is easy to check that
 \beqnn
\xi_n(g_{\lambda,z}) := \int_{\mbb{R}_+\times \mbb{R}} g_{\lambda,z}
d\xi_n
 =
\lambda_1a_n^{-2} S_{1,n}(a_nz,\infty) + \lambda_2c_n^{-1} S_{2,n}(c_n
z,\infty).
 \eeqnn
On the other hand, one can see the mapping from $M(\bar{\mbb{R}}^2_0)$
into $\mbb{R}$ defined by
 \beqnn
N := \sum_{k=1}^\infty\delta_{(x_{1,k},x_{2,k})}
 \mapsto
N(g_{\lambda,z}) := \int_{\mbb{R}_+\times \mbb{R}} g_{\lambda,z}dN
 \eeqnn
is a.s.\ continuous with respect to distribution of the limit process
$\xi$ in Theorem~\ref{t4.9}. By the continuous mapping theorem, as $n\to
\infty$, we have $\xi_n(g_{\lambda,z})\overset{\rm d} {\longrightarrow}
\xi(g_{\lambda,z})$, and hence
 \beqnn
\mbf{E}[\exp(i\xi_n(g_{\lambda,z}))]
 \to
\mbf{E}[\exp(i\xi(g_{\lambda,z}))],
 \eeqnn
where the right hand side is equal to
 \beqnn
\ar\ar\exp\bigg\{-\frac{a\sigma^\alpha}{\alpha^2b^2\Gamma(-\alpha)}
\int_0^\infty \mbf{E}\Big[\exp\Big\{i\lambda_1y^2 \sum_{j=2}^\infty
e^{-2b(j-1)}1_{\{ye^{-b(j-1)}>z\}}\Big\} \cr
 \ar\ar\qquad\qquad
\exp\Big\{i\lambda_2y^{(\alpha+1)/\alpha} \sum_{j=2}^\infty
e^{-b(j-1)(\alpha+1)/\alpha}V_j1_{\{|e^{-b(j-1)
(\alpha+1)/\alpha}V_j|>z\}}\Big\}\Big] \cr
 \ar\ar\qquad\qquad
\mbf{E}\Big(1-\exp\{i\lambda_1y^2+i\lambda_2y^{(\alpha+1)/\alpha}
V_1\}\Big) \frac{dy}{y^{\alpha+1}}\bigg\}.
 \eeqnn
Then we can use dominated convergence theorem to see that, as
$z\rightarrow 0$,
 \beqnn
\ar\ar\exp\bigg\{-\frac{a\sigma^\alpha}{\alpha^2b^2\Gamma(-\alpha)}
\int_0^\infty \mbf{E}\Big(1-\exp\{i\lambda_1y^2 + i\lambda_2
y^{(\alpha+1)/\alpha} V_1\}\Big) \cr
 \ar\ar\quad
\mbf{E}\Big[\exp\Big\{i \sum_{j=2}^\infty \big(\lambda_1y^2 e^{-2b(j-1)} +
\lambda_2 y^{(\alpha+1)/\alpha} e^{-b(j-1)(\alpha+1)/\alpha}
V_j\big)\Big\}\Big] \frac{dy}{y^{\alpha+1}}\bigg\}.
 \eeqnn
Since the sequence $\{V_1,V_2,\cdots\}$ is i.i.d., the above quantity is
equal to (\ref{4.11}). Note that $\mbf{E}(V_1) = 0$. Then the theorem
follows by Lemmas~\ref{t4.10} and \ref{t4.11}. \qed


\section{Asymptotics of the estimators}

\setcounter{equation}{0}

In this section, we investigate the asymptotics of the estimators for the
SCIR-model. The results are presented in a number of theorems. In fact,
we shall first study the asymptotics of the estimators of the parameters
$(\gamma,\rho)$ defined in (\ref{1.4}). Their CLSEs can be obtained by
minimizing the sum of squares in (\ref{1.8}). They are given by
 \beqlb\label{5.1}
\hat{\gamma}_n = \frac {\sum_{k=1}^nX_{k-1}\sum_{k=1}^nX_k -
n\sum_{k=1}^nX_{k-1}X_k} {\big(\sum_{k=1}^nX_{k-1}\big)^2 -
n\sum_{k=1}^nX_{k-1}^2}
 \eeqlb
and
 \beqlb\label{5.2}
\hat{\rho}_n
 \ar=\ar
\frac{1}{n}\Big[\sum_{k=1}^nX_k - \hat{\gamma}_n \sum_{k=1}^n
X_{k-1}\Big].
 \eeqlb
By minimizing the weighted sum in (\ref{1.11}), we obtain the WCLSEs of
the parameters:
 \beqlb\label{5.3}
\check{\gamma}_n
 =
\frac{\sum_{k=1}^nX_k\sum_{k=1}^n\frac{1}{X_{k-1}+1}-n \sum_{k=1}^n
\frac{X_k}{X_{k-1}+1}} {\sum_{k=1}^n(X_{k-1}+1)\sum_{k=1}^n
\frac{1}{X_{k-1}+1}-n^2}.
 \eeqlb
and
 \beqlb\label{5.4}
\check{\rho}_n
 \ar=\ar
\frac{1}{n}\Big[\sum_{k=1}^nX_k - \check{\gamma}_n \sum_{k=1}^n
X_{k-1}\Big].
 \eeqlb
In view of Proposition~\ref{t2.3} and the above expressions, in the
discussions of the above estimators it suffices to consider a stationary
realization $\{X_t: t\ge 0\}$ of the SCIR-model.

\blemma\label{t5.1} We have, as $n\to \infty$,
 \beqlb\label{5.5}
\frac{1}{n}\sum_{k=1}^nX_k\overset{\rm a.s.}{\longrightarrow} \frac{a}{b},
 \quad
\frac{1}{n}\sum_{k=1}^n\frac{1}{1+X_{k-1}} \overset{\rm
a.s.}{\longrightarrow} \lambda,
 \eeqlb
and
 \beqlb\label{5.6}
\frac{1}{n}\sum_{k=1}^n\frac{X_{k}}{X_{k-1}+1}
 \overset{\rm a.s.}{\longrightarrow}
\rho\lambda + \gamma(1-\lambda)
 \eeqlb
where
 \beqnn
\lambda = \mbf{E}\Big(\frac{1}{1+X_0}\Big).
 \eeqnn
 \elemma

\proof By Theorem~\ref{t2.6}, the process $\{X_t\}$ is exponentially
ergodic and thus strong mixing, so the tail $\sigma$-algebra of the
process is trivial; see, e.g., Durrett (1996, p.351). Recall that
$\mbf{E}(X_0) = a/b$. In view of (\ref{1.6}), we have
 \beqnn
\mbf{E}\Big(\frac{X_1}{1+X_0}\Big)
 =
\rho\mbf{E}\Big(\frac{1}{1+X_0}\Big) +
\gamma\mbf{E}\Big(\frac{X_0}{1+X_0}\Big)
 =
\rho\lambda + \gamma(1-\lambda).
 \eeqnn
Then the result follows by Birkhoff's ergodic theorem; see, e.g., Durrett
(1996, p.341). \qed

\ttheorem\label{t5.2} The estimators $(\check{\rho}_n, \check{\gamma}_n)$
are strongly consistent and, as $n\to \infty$, $n^{(\alpha-1)/\alpha}
\big(\check{\gamma}_n-\gamma, \check{\rho}_n-\rho\big)$ converges in
distribution to
 \beqnn
F^{-1} (U_1, U_2)
 \Big(\begin{array}{cc}
\lambda & \lambda-1 \ccr -1 & ab^{-1}
 \end{array}\Big)
 =
F^{-1}\big(\lambda U_1 - U_2, (\lambda-1)U_1 + ab^{-1}U_2\big),
 \eeqnn
where $F=(1+ab^{-1})\lambda-1$ and $(U_1,U_2)$ is an $\alpha$-stable
random vector with characteristic function given by (\ref{4.7}).
\etheorem

\proof We first remark that (\ref{4.7}) defines a Gaussian random vector
$(U_1,U_2)$ when $\alpha=2$. In view of (\ref{5.3}) and (\ref{5.4}), the
results of Lemma~\ref{t5.1} imply that
 \beqnn
\check{\gamma}_n
 \ar\overset{\rm a.s.}{\longrightarrow}\ar
\frac{ab^{-1}\lambda - \rho\lambda - \gamma(1-\lambda)}
{(1+ab^{-1})\lambda - 1}
 =
\frac{ab^{-1}\gamma\lambda - \gamma(1-\lambda)} {(1+ab^{-1})\lambda - 1}
= \gamma
 \eeqnn
and
 \beqnn
\check{\rho}_n
 \ar\overset{\rm a.s.}{\longrightarrow}\ar
\frac{a}{b} (1-\gamma) = \rho.
 \eeqnn
Those give the strong consistency of $\check{\rho}_n$ and
$\check{\gamma}_n$. Write
 \beqlb\label{5.7}
\gamma
 \ar=\ar
\frac{\sum_{k=1}^n(\gamma X_{k-1}+\gamma)\sum_{k=1}^n\frac{1}{X_{k-1}+1}
- n\sum_{k=1}^n\frac{\gamma X_{k-1}+\gamma}{X_{k-1}+1}}
{\sum_{k=1}^n(X_{k-1}+1)\sum_{k=1}^n\frac{1}{X_{k-1}+1}-n^2}\cr
 \ar=\ar
\frac{\sum_{k=1}^n\gamma X_{k-1}\sum_{k=1}^n\frac{1}{X_{k-1}+1} -
n\sum_{k=1}^n\frac{\gamma X_{k-1}}{X_{k-1}+1}}
{\sum_{k=1}^n(X_{k-1}+1)\sum_{k=1}^n\frac{1}{X_{k-1}+1}-n^2}\cr
 \ar=\ar
\frac{\sum_{k=1}^n(\gamma X_{k-1}+\rho)\sum_{k=1}^n\frac{1}{X_{k-1}+1} -
n\sum_{k=1}^n\frac{\gamma X_{k-1}+\rho}{X_{k-1}+1}}
{\sum_{k=1}^n(X_{k-1}+1)\sum_{k=1}^n\frac{1}{X_{k-1}+1}-n^2}\cr
 \ar=\ar
\frac{\sum_{k=1}^n(X_k-\varepsilon_k)\sum_{k=1}^n\frac{1}{X_{k-1}+1} -
n\sum_{k=1}^n\frac{X_k-\varepsilon_k}{X_{k-1}+1}} {\sum_{k=1}^n
(X_{k-1}+1)\sum_{k=1}^n\frac{1}{X_{k-1}+1}-n^2},
 \eeqlb
where the last equality follows from (\ref{1.6}). Combining this with
(\ref{5.3}) we have
 \beqlb\label{5.8}
\check{\gamma}_n-\gamma
 \ar=\ar
\frac{\sum_{k=1}^n\frac{1}{X_{k-1}+1}\sum_{k=1}^n\varepsilon_k
-n\sum_{k=1}^n\frac{\varepsilon_k}{X_{k-1}+1}}
{\sum_{k=1}^n(X_{k-1}+1)\sum_{k=1}^n\frac{1}{X_{k-1}+1}-n^2}.
 \eeqlb
By (\ref{1.6}) and (\ref{5.4}) it is easy to see that
 \beqnn
\check{\rho}_n
 \ar=\ar
\frac{\sum_{k=1}^n(X_k-\gamma X_{k-1})-(\check{\gamma}_n-\gamma)
\sum_{k=1}^nX_{k-1}}{n}\cr
 \ar=\ar
\frac{\sum_{k=1}^n(\rho+\varepsilon_k) - (\check{\gamma}_n - \gamma)
\sum_{k=1}^nX_{k-1}} {n}.
 \eeqnn
Using (\ref{5.8}) we have
 \beqlb\label{5.9}
\check{\rho}_n-\rho
 \ar=\ar
\frac{\sum_{k=1}^n\varepsilon_k-(\check{\gamma}_n-\gamma)
\sum_{k=1}^nX_{k-1}} {n}\cr
 \ar=\ar
\frac{\sum_{k=1}^n\frac{\varepsilon_k}{X_{k-1}+1}\sum_{k=1}^nX_{k-1} -
\sum_{k=1}^n\varepsilon_k\sum_{k=1}^n\frac{X_{k-1}}{X_{k-1}+1}}
{\sum_{k=1}^n(X_{k-1}+1)\sum_{k=1}^n\frac{1}{X_{k-1}+1}-n^2}.
 \eeqlb
We can rewrite (\ref{5.8}) and (\ref{5.9}) into the matrix form as
 \beqnn
n^{(\alpha-1)/\alpha}\big(\check{\gamma}_n - \gamma, \check{\rho}_n -
\rho\big)
 =
F_n^{-1}\mbf{U}_n\mbf{B}_n,
 \eeqnn
where
 \beqnn
F_n = \frac{1}{n^2}\sum_{k=1}^n(X_{k-1}+1)
\sum_{k=1}^n\frac{1}{X_{k-1}+1} - 1, \eeqnn
 \beqnn
\mbf{B}_n = \frac{1}{n} \left(\begin{array}{cc}
\sum_{k=1}^n\frac{1}{X_{k-1}+1} & -\sum_{k=1}^n\frac{X_{k-1}}{X_{k-1}+1}
 \ccr
-n & \sum_{k=1}^n X_{k-1}
\end{array}\right),
 \eeqnn
and
 \beqnn
\mbf{U}_n = \frac{1}{n^{1/\alpha}} \bigg(\sum_{k=1}^n\varepsilon_k,
\sum_{k=1}^n \frac{\varepsilon_k}{X_{k-1}+1} \bigg).
 \eeqnn
From (\ref{5.5}) it follows that $F_n\overset{\rm a.s.}{\longrightarrow}
F$ and
 \beqnn
\mbf{B}_n\overset{\rm a.s.}{\longrightarrow} \mbf{B} :=
 \Big(\begin{array}{cc}
\lambda & \lambda-1
 \ccr
-1 & ab^{-1}
\end{array}\Big).
 \eeqnn
In the case $1<\alpha<2$, we have $\mbf{U}_n\overset{\rm
d}{\longrightarrow} (U_1,U_2)$ by Theorem~\ref{t4.5}. In the case
$\alpha=2$, we have
 \beqnn
\mbf{E}[\varepsilon_1^2]
 =
\sigma^2\mbf{E}\Big[\int_0^1 e^{-2b(1-s)} X_sds\Big]
 =
\frac{a\sigma^2}{b}\int_0^1 e^{-2b(1-s)} ds
 <
\frac{a\sigma^2}{b}.
 \eeqnn
Then, for any $\delta>0$,
 \beqnn
\ar\ar\frac{1}{n}\sum_{k=1}^n \mbf{E}\Big[\Big(\lambda_1 +
\frac{\lambda_2} {X_{k-1}+1}\Big)^2\varepsilon_k^2
1_{\big\{\big(\lambda_1 + \frac{\lambda_2}
{X_{k-1}+1}\big)|\varepsilon_k|> \delta\sqrt{n}\big\}}\Big] \cr
 \ar\ar\qquad\qquad
= \mbf{E}\Big[\Big(\lambda_1 + \frac{\lambda_2} {X_0+1}\Big)^2
\varepsilon_1^2 1_{\big\{\big(\lambda_1 + \frac{\lambda_2}
{X_0+1}\big)|\varepsilon_1|> \delta\sqrt{n}\big\}}\Big]
 \eeqnn
tends to zero as $n\to \infty$. Furthermore, by the ergodicity theorem,
as $n\to \infty$,
 \beqlb\label{5.10}
\ar\ar\frac{1}{n}\sum_{k=1}^n \mbf{E}\Big[\Big(\lambda_1 +
\frac{\lambda_2} {X_{k-1}+1}\Big)^2 \varepsilon_k^2\Big|
\mcr{F}_{k-1}\Big] \cr
 \ar\ar\qquad\qquad\overset{\rm a.s.}\longrightarrow
\sigma^2\mbf{E}\Big[\Big(\lambda_1 + \frac{\lambda_2} {1+X_0}\Big)^2
\big(p_2 + q_2X_0\big)\Big].
 \eeqlb
Then a martingale convergence theorem implies that, as $n\to \infty$,
 \beqnn
\frac{1}{\sqrt{n}}\sum_{k=1}^n\Big(\lambda_1 + \frac{\lambda_2}
{X_{k-1}+1}\Big)\varepsilon_k
 \eeqnn
converges in distribution to a Gaussian random variable with mean zero
and variance given by the right-hand side of (\ref{5.10}); see, e.g.,
Durrett (1996, p.417). It follows that
 \beqnn
\mbf{E}\Big[\exp\Big\{\frac{i}{\sqrt{n}}\sum_{k=1}^n\Big(\lambda_1 +
\frac{\lambda_2} {X_{k-1}+1}\Big)\varepsilon_k\Big\}\Big]
 \eeqnn
converges to the right-hand side of (\ref{4.7}) with $\alpha=2$. Since
$(\lambda_1,\lambda_2)\in \mbb{R}^2$ can be arbitrary in the above, we
also conclude $\mbf{U}_n\overset{\rm d}{\longrightarrow} (U_1,U_2)$. That
proves the desired convergence. \qed

\ttheorem\label{t5.3} The estimators $(\check{b}_n, \check{a}_n)$ are
strongly consistent and as $n\to \infty$, $n^{(\alpha-1)/\alpha}
(\check{b}_n - b, \check{a}_n-a)$ converges in distribution to
 \beqnn
F^{-1}\big(e^{b}(U_2-\lambda U_1), (1-e^{-b})^{-1} [a\lambda +
b(\lambda-1)] U_1 + ab^{-1}e^b(U_2-\lambda U_1)\big),
 \eeqnn
where $F = (1+ab^{-1})\lambda-1$ and $(U_1,U_2)$ is an $\alpha$-stable
random vector with characteristic function given by (\ref{4.7}).
\etheorem

\proof The strong consistency of $\check{b}_n$ and $\check{a}_n$ follows
from that of $\check{\rho}_n$ and $\check{\gamma}_n$. By the relations in
(\ref{1.4}), we have, as $n\to \infty$,
 \beqlb\label{5.11}
(\check{\gamma}_n - \gamma)
 =
e^{-\check{b}_n} - e^{-b}
 =
-(\check{b}_n - b)e^{-b} + o(\check{b}_n - b)
 \eeqlb
and
 \beqlb\label{5.12}
\check{a}_n-a
 \ar=\ar
\frac{\check{\rho}_n\check{b}_n}{1-e^{-\check{b}_n}} - \frac{\rho
b}{1-e^{-b}}
 =
\frac{\check{\rho}_n\check{b}_n(1-e^{-b}) - \rho b(1-e^{-\check{b}_n})}
{(1-e^{-\check{b}_n}) (1-e^{-b})} \cr
 \ar=\ar
\frac{\check{b}_n(\check{\rho}_n-\rho)} {1-e^{-\check{b}_n}} +
\frac{\rho(\check{b}_n-b)} {1-e^{-\check{b}_n}} + \frac{\rho
b(e^{-\check{b}_n} - e^{-b})} {(1-e^{-\check{b}_n})(1-e^{-b})} \cr
 \ar=\ar
\frac{b(\check{\rho}_n-\rho)} {1-e^{-b}} - \frac{ae^b(\check{\gamma}_n -
\gamma)} {b} + \frac{a(\check{\gamma}_n - \gamma)} {1-e^{-b}} +
o(\check{b}_n - b).
 \eeqlb
Then the desired convergence follows from Theorem~\ref{t5.2}. \qed

\ttheorem\label{t5.4} The estimators $(\hat{\rho}_n, \hat{\gamma}_n)$ are
weakly consistent. Moreover, if $1<\alpha< (1+\sqrt{5})/2$, then, as
$n\to \infty$,
 \beqlb\label{5.13}
n^{(\alpha-1)/\alpha^2}(\hat{\gamma}_n - \gamma, \hat{\rho}_n-\rho)
 \overset{\rm d}{\longrightarrow}
S_1^{-1}S_2(1, -ab^{-1}),
 \eeqlb
where $(S_1,S_2)$ has characteristic function given by (\ref{4.11}).
 \etheorem

\proof By (\ref{1.6}) and (\ref{5.1}) we have
 \beqnn
\hat{\gamma}_n - \gamma
 \ar=\ar
\frac{\sum_{k=1}^nX_{k-1}\sum_{k=1}^n(\varepsilon_k+\rho) - n\sum_{k=1}^n
X_{k-1}(\varepsilon_k+\rho)} {\big(\sum_{k=1}^nX_{k-1}\big)^2 -
n\sum_{k=1}^nX_{k-1}^2} \cr
 \ar=\ar
\frac{\sum_{k=1}^nX_{k-1}\sum_{k=1}^n\varepsilon_k - n\sum_{k=1}^n
X_{k-1}\varepsilon_k} {\big(\sum_{k=1}^nX_{k-1}\big)^2 -
n\sum_{k=1}^nX_{k-1}^2}.
 \eeqnn
Then using (\ref{5.2}) we get
 \beqnn
\hat{\rho}_n - \rho
 \ar=\ar
\frac{\sum_{k=1}^nX_k - \hat{\gamma}_n\sum_{k=1}^n X_{k-1}}{n} - \rho \cr
 \ar=\ar
\frac{\sum_{k=1}^n(X_k-\gamma X_{k-1}) - (\hat{\gamma}_n-\gamma)
\sum_{k=1}^n X_{k-1}}{n} - \rho \cr
 \ar=\ar
\frac{\sum_{k=1}^n(\varepsilon_k+\rho) - (\hat{\gamma}_n-\gamma)
\sum_{k=1}^n X_{k-1}}{n} - \rho \cr
 \ar=\ar
\frac{\sum_{k=1}^n \varepsilon_k - (\hat{\gamma}_n-\gamma) \sum_{k=1}^n
X_{k-1}}{n} \cr
 \ar=\ar
\frac{\sum_{k=1}^nX_{k-1}\sum_{k=1}^nX_{k-1}\varepsilon_k - \sum_{k=1}^n
X_{k-1}^2\sum_{k=1}^n\varepsilon_k} {\big(\sum_{k=1}^nX_{k-1}\big)^2 -
n\sum_{k=1}^nX_{k-1}^2}.
 \eeqnn
By Theorem~\ref{t4.12} it is easy to see that
 \beqnn
\frac{\sum_{k=1}^nX_{k-1}\varepsilon_k}
{\sum_{k=1}^nX_{k-1}^2}\overset{\rm p}\longrightarrow 0.
 \eeqnn
Then we have $\hat{\rho}_n-\rho\overset{\rm p}\longrightarrow 0$ and
$\hat{\gamma}_n-\gamma\overset{\rm p}\longrightarrow 0$, giving the weak
consistency of $(\hat{\rho}_n, \hat{\gamma}_n)$. From the above relations
it follows that
 {\small\beqnn
n^{(\alpha-1)/\alpha^2}(\hat{\gamma}_n - \gamma)
 =
\frac{\frac{1}{n^{1+(\alpha+1)/\alpha^2}}\sum_{k=1}^nX_{k-1}
\sum_{k=1}^n\varepsilon_k - \frac{1}{n^{(\alpha+1)/\alpha^2}}
\sum_{k=1}^nX_{k-1}\varepsilon_k} {\frac{1}{n^{1+2/\alpha}}
\big(\sum_{k=1}^nX_{k-1}\big)^2-\frac{1}{n^{2/\alpha}}\sum_{k=1}^n
X_{k-1}^2}
 \eeqnn}\!\!
and
 {\small\beqnn
n^{(\alpha-1)/\alpha^2}(\hat{\rho}_n - \rho)
 =
\frac{\frac{1}{n^{1+(\alpha+1)/\alpha^2}}\big(\sum_{k=1}^nX_{k-1}
\sum_{k=1}^nX_{k-1}\varepsilon_k - \sum_{k=1}^nX_{k-1}^2
\sum_{k=1}^n\varepsilon_k\big)} {\frac{1}{n^{1+2/\alpha}}
\big(\sum_{k=1}^nX_{k-1}\big)^2 - \frac{1}{n^{2/\alpha}}
\sum_{k=1}^nX_{k-1}^2}.
 \eeqnn}\!\!
Take any constant $0<\delta< [1\land (\alpha-1)^2]/\alpha^2$. We can
rewrite the above relations into the matrix form
 \beqnn
n^{(\alpha-1)/\alpha}\big(\hat{\gamma}_n - \gamma, \hat{\rho}_n -
\rho\big)
 =
T_n^{-1}\mbf{S}_n\mbf{A}_n,
 \eeqnn
where
 \beqnn
T_n = \frac{1}{n^{1+2/\alpha}} \Big(\sum_{k=1}^nX_{k-1}\Big)^2 -
\frac{1}{n^{2/\alpha}} \sum_{k=1}^nX_{k-1}^2,
 \eeqnn
 \beqnn
\mbf{A}_n = \bigg(\begin{array}{cc} \frac{1}{n^{1+1/\alpha^2-\delta}}
\sum_{k=1}^n X_{k-1} & - \frac{1}{n^{1+1/\alpha^2-\delta}} \sum_{k=1}^n
X_{k-1}^2
 \ccr
-1 & \frac{1}{n}\sum_{k=1}^n X_{k-1}
\end{array}\bigg),
 \eeqnn
and
 \beqnn
\mbf{S}_n = \bigg(\frac{1}{n^{1/\alpha+\delta}}
\sum_{k=1}^n\varepsilon_k, \frac{1}{n^{(\alpha+1)/\alpha^2}}\sum_{k=1}^n
X_{k-1}\varepsilon_k \bigg).
 \eeqnn
If $\alpha^2< \alpha+1$, by (\ref{5.5}) and Theorem~\ref{t4.12}, we have
$T_n\overset{\rm d}{\longrightarrow} -S_1$ and
 \beqnn
\mbf{A}_n\overset{\rm d}{\longrightarrow} \mbf{A} :=
 \bigg(\begin{array}{cc}
0 & 0
 \ccr
-1 & ab^{-1}
\end{array}\bigg).
\eeqnn By Theorems~\ref{t4.5} and~\ref{t4.12} we have
$\mbf{S}_n\overset{\rm d}{\longrightarrow} (0,S_2)$. Then (\ref{5.13})
holds. \qed

\ttheorem\label{t5.5} The estimators $(\hat{b}_n, \hat{a}_n)$ are weakly
consistent. Moreover, if $1<\alpha< (1+\sqrt{5})/2$, then, as $n\to
\infty$,
 \beqnn
n^{(\alpha-1)/\alpha^2}(\hat{b}_n-b, \hat{a}_n-a)
 \overset{\rm d}{\longrightarrow}
- e^{b}(1, ab^{-1})S_1^{-1}S_2,
 \eeqnn
where $(S_1,S_2)$ the characteristic function given by (\ref{4.11}).
 \etheorem

\proof The weak consistency of $\hat{b}_n$ and $\hat{a}_n$ follows from
that of $\hat{\rho}_n$ and $\hat{\gamma}_n$. The relations (\ref{5.11})
and (\ref{5.12}) still hold when the ``checks'' are replaced by ``hats''.
Then the desired result follows from Theorem~\ref{t5.4}. \qed

\ttheorem\label{t5.6} The estimator $\hat{\sigma}_n$ for $\sigma$ defined
by (\ref{1.14}) is weakly consistent. \etheorem

\proof For $p\in (0,\alpha)$ and $\delta\in (0, \min\{1-1/\alpha,
1/\alpha^2\})$, we need to show that, as $n\to \infty$,
 \beqlb\label{5.14}
\frac{1}{n^{1-p/\alpha}}\sum_{k=1}^n\Big|\frac{X_{t_k}-X_{t_{k-1}}}
{X_{t_{k-1}}^{1/\alpha} + n^{-\delta}}\Big|^p
 \overset{\rm p}{\longrightarrow}
\sigma^p\mbf{E}[|Z_1|^p].
 \eeqlb
By equation (\ref{1.2}), we have
 \beqnn
\mbox{l.h.s. of (\ref{5.14})}
 \ar=\ar
\frac{1}{n^{1-p/\alpha}}\sum_{k=1}^n\Big(\Big|\frac{X_{t_k}-X_{t_{k-1}}}
{X_{t_{k-1}}^{1/\alpha}+n^{-\delta}}\Big|^p -
\Big|\frac{\int_{t_{k-1}}^{t_k}\sigma X_{s-}^{1/\alpha}dZ_s}
{X_{t_{k-1}}^{1/\alpha}+n^{-\delta}}\Big|^p\Big) \cr
 \ar\ar
+ \frac{1}{n^{1-p/\alpha}}\sum_{k=1}^n\Big(\Big|\int_{t_{k-1}}^{t_k}
\frac{\sigma X_{s-}^{1/\alpha}} {X_{t_{k-1}}^{1/\alpha}+n^{-\delta}}
dZ_s\Big|^p \cr
 \ar\ar\qquad\qquad\qquad\qquad
-\, \Big|\int_{t_{k-1}}^{t_k}\frac{\sigma X_{t_{k-1}-}^{1/\alpha}}
{X_{t_{k-1}}^{1/\alpha} + n^{-\delta}}dZ_s\Big|^p\Big) \cr
 \ar\ar
+ \frac{1}{n^{1-p/\alpha}}\sum_{k=1}^n\Big(\Big|\int_{t_{k-1}}^{t_k}
\frac{\sigma X_{t_{k-1}-}^{1/\alpha}} {X_{t_{k-1}}^{1/\alpha}+n^{-\delta}}
dZ_s\Big|^p - \Big|\int_{t_{k-1}}^{t_k} \sigma dZ_s\Big|^p\Big) \cr
 \ar\ar
+ \frac{1}{n^{1-p/\alpha}}\sum_{k=1}^n \Big|\int_{t_{k-1}}^{t_k} \sigma
dZ_s\Big|^p.
 \eeqnn
It is easy to see that
 \beqnn
J_1 \ar:=\ar \frac{1}{n^{1-p/\alpha}}\sum_{k=1}^n
\bigg|\Big|\frac{X_{t_k}-X_{t_{k-1}}} {X_{t_{k-1}}^{1/\alpha} +
n^{-\delta}}\Big|^p - \Big|\frac{\int_{t_{k-1}}^{t_k} \sigma
X_{s-}^{1/\alpha}dZ_s} {X_{t_{k-1}}^{1/\alpha} + n^{-\delta}}\Big|^p\bigg|
\cr
 \ar\le\ar
\frac{1}{n^{1-p/\alpha}}\sum_{k=1}^n\Big|\frac{-b\int_{t_{k-1}}^{t_k}
X_sds+a(t_k-t_{k-1})} {X_{t_{k-1}}^{1/\alpha}+n^{-\delta}}\Big|^p \cr
 \ar\le\ar
n^{p/\alpha}\Big(b\sup_{s\in[0,1]}X_s + a\Big)^pn^{p(\delta-1)},
 \eeqnn
which goes to zero a.s.\ as $n\to \infty$. Observe that
 \beqnn
J_2 \ar:=\ar \frac{1}{n^{1-p/\alpha}}\sum_{k=1}^n \bigg|\Big|
\int_{t_{k-1}}^{t_k}\frac{\sigma X_{s-}^{1/\alpha}}
{X_{t_{k-1}}^{1/\alpha}+n^{-\delta}}dZ_s\Big|^p -
\Big|\int_{t_{k-1}}^{t_k}\frac{\sigma X_{t_{k-1}-}^{1/\alpha}}
{X_{t_{k-1}}^{1/\alpha} +n^{-\delta}}dZ_s\Big|^p\bigg| \cr
 \ar\le\ar
\frac{1}{n^{1-p/\alpha-p\delta}}\sum_{k=1}^n
\Big|\int_{t_{k-1}}^{t_k}\sigma(X_{s-}^{1/\alpha} -
X_{t_{k-1}}^{1/\alpha})dZ_s\Big|^p.
 \eeqnn
By equation (\ref{1.2}) and Lemma~\ref{t2.8},
 \beqnn
\ar\ar\mbf{E}\Big[\sup_{s\in[t_{k-1},t_k]}|X_s-X_{t_{k-1}}|\Big] \cr
 \ar\ar\qquad
\le \mbf{E}\Big[b\int_{t_{k-1}}^{t_k}X_sds\Big] + \frac{a}{n} +
\mbf{E}\Big[\sup_{s\in[t_{k-1},t_k]}\Big|\int_{t_{k-1}}^s\sigma
X_{u-}^{1/\alpha}dZ_u\Big|\Big]\cr
 \ar\ar\qquad
\le \Big(a+b\mbf{E}\Big[\sup_{s\in[0,1]}X_s\Big]\Big)\frac{1}{n} +
C_1\sigma \mbf{E}\Big[\Big(\int_{t_{k-1}}^{t_k} X_sds\Big)^{1/\alpha}\Big]
\cr
 \ar\ar\qquad
\le C_2\Big(\frac{1}{n}+\frac{1}{n^{1/\alpha}}\Big).
 \eeqnn
Then by Lemma~\ref{t2.8},
 \beqnn
\ar\ar \frac{1}{n^{1-p\delta-p/\alpha}} \mbf{E}\Big[\sum_{k=1}^n
\Big|\int_{t_{k-1}}^{t_k}\sigma(X_{s-}^{1/\alpha}-
 X_{t_{k-1}}^{1/\alpha})dZ_s\Big|^p \Big]\cr
 \ar\ar\qquad
\le \frac{L}{n^{1-p\delta-p/\alpha}}\sum_{k=1}^n
\mbf{E}\Big[\Big(\int_{t_{k-1}}^{t_k}
|X_s-X_{t_{k-1}}|ds\Big)^{p/\alpha}\Big]\cr
 \ar\ar\qquad
\le \frac{L}{n^{1-p\delta-p/\alpha}}\sum_{k=1}^n
\Big(\int_{t_{k-1}}^{t_k} \mbf{E}[|X_s-X_{t_{k-1}}|]ds\Big)^{p/\alpha}
\cr
 \ar\ar\qquad
\le C_3n^{p\delta-p/\alpha^2}.
 \eeqnn
The right-hand side tends to zero as $n\to \infty$. It follows that
$J_2\overset{\rm p}{\longrightarrow} 0$ as $n\to \infty$. Using the
stationarity of $\{X_t\}$ and the self-similarity of $\{Z_t\}$ we have
 \beqnn
\frac{1}{n^{1-p/\alpha}}\mbf{E}\Big[\sum_{k=1}^n\Big|\int_{t_{k-1}}^{t_k}
\frac{n^{-\delta}}{X_{t_{k-1}}^{1/\alpha}+n^{-\delta}} dZ_s\Big|^p\Big]
 =
\mbf{E}\Big[\Big(\frac{n^{-\delta}}{X_0^{1/\alpha} +
n^{-\delta}}\Big)^p\Big]\mbf{E}[|Z_1|^p],
 \eeqnn
which goes to zero by dominated convergence theorem. Observe that
 \beqnn
J_3 \ar:=\ar \frac{1}{n^{1-p/\alpha}}\sum_{k=1}^n\bigg|\Big|
\int_{t_{k-1}}^{t_k}\frac{\sigma X_{t_{k-1}}^{1/\alpha}}
{X_{t_{k-1}}^{1/\alpha}+n^{-\delta}}dZ_s\Big|^p -
\Big|\int_{t_{k-1}}^{t_k}\sigma dZ_s\Big|^p\bigg| \cr
 \ar\le\ar
\frac{L}{n^{1-p/\alpha}}\sum_{k=1}^n\big|\int_{t_{k-1}}^{t_k}
\frac{n^{-\delta}}{X_{t_{k-1}}^{1/\alpha}+n^{-\delta}}dZ_s\big|^p.
 \eeqnn
Then $J_3\overset{\rm p}{\longrightarrow} 0$ as $n\to \infty$. On the
other hand, it is not hard to see that
 \beqnn
\frac{1}{n^{1-p/\alpha}}\sum_{k=1}^n{\sigma}^{p}|Z_{t_k}-Z_{t_{k-1}}|^p
 \overset{\rm d}{=}
\frac{1}{n}\sum_{k=1}^n{\sigma}^{p}|Z_k-Z_{k-1}|^p
 \overset{\rm p}{\longrightarrow}
\sigma^{p}\mbf{E}[|Z_1|^p].
 \eeqnn
Thus we have (\ref{5.14}). \qed

\noindent\textbf{Acknowledgements.} We would like to thank Professors
Thomas Mikosch and Gennady Samorodnitsky for enlightening discussions on
regular variations and related properties of stochastic processes. We are
grateful to Professor Matyas Barczy for his careful reading of earlier
versions of this paper and pointing out many typos and errors.

\bigskip


\medskip

\noindent School of Mathematical Sciences \\
Beijing Normal University \\
Beijing 100875, P.\,R.\ China \\
E-mail: {\tt lizh@bnu.edu.cn} \\
~

\smallskip

\noindent School of Mathematical Sciences and LPMC \\
Nankai University \\
Tianjin 300071, P.\,R.\ China \\
E-mail: {\tt mach@nankai.edu.cn} \\
~

\end{document}